\documentclass[12pt]{article}

\usepackage{latexsym}

\begin{document}


\newcommand{\deu}[1]{}
\newcommand{\eng}[1]{#1}

%

\newcommand{\EE}{\mathop{\rm I\! E}\nolimits}
\newcommand{\E}{\mathop{\rm E}\nolimits}
\newcommand{\I}{\mathop{\rm Im  }\nolimits}
\newcommand{\R}{\mathop{\rm Re  }\nolimits}
\newcommand{\CC}{\mathop{\rm C\!\!\! I}\nolimits}
\newcommand{\FF}{\mathop{\rm I\! F}\nolimits}
\newcommand{\KK}{\mathop{\rm I\! K}\nolimits}
\newcommand{\LL}{\mathop{\rm I\! L}\nolimits}
\newcommand{\MM}{\mathop{\rm I\! M}\nolimits}
\newcommand{\NN}{\mathop{\rm I\! N}\nolimits}
\newcommand{\PP}{\mathop{\rm I\! P}\nolimits}
\newcommand{\QQ}{\mathop{\rm I\! Q}\nolimits}
\newcommand{\RR}{\mathop{\rm I\! R}\nolimits}
\newcommand{\ZZ}{\mathop{\rm Z\!\!Z}\nolimits}
\newcommand{\integer}{\mathop{\rm int}\nolimits}
\newcommand{\erf}{\mathop{\rm erf}\nolimits}
\newcommand{\diag}{\mathop{\rm diag}\nolimits}
\newcommand{\fl}{\mathop{\rm fl}\nolimits}
\newcommand{\eps}{\mathop{\rm eps}\nolimits}
\newcommand{\var}{\mathop{\rm var}\nolimits}


\newcommand{\bild}[3]{{        
  \unitlength1mm
  \begin{figure}[ht]
  \begin{picture}(120,#1)\end{picture}
  \caption{\label{#3}#2}
  \end{figure}
}}

\newcommand{\dontdothat}[1]{}

\newcommand{\ring}{{\cal K}}

\newcommand{\pfeil}{\rightarrow}

\newcommand{\HA}{{\rm HA}}
\newcommand{\HB}{{\rm HB}}
\newcommand{\BA}{{\rm BA}}
\newcommand{\ZA}{{\rm ZA}}
\newcommand{\ZB}{{\rm ZB}}
\newcommand{\HBB}{{\rm\cal HB}}
\newcommand{\BAB}{{\rm\cal BA}}
\newcommand{\BBfull}{{\rm BB}}
\newcommand{\BB}{{\rm BB}}
\newcommand{\BBB}{{\rm\cal BB}}
\newcommand{\BBs}{{\rm {\widetilde{BB}}}}
\newcommand{\tr}{{\rm tr}}
\newcommand{\Tr}{{\rm Tr}}
\newcommand{\iso}{\stackrel{\sim}{=}}

\newcommand{\kat}{{\cal C}}
\newcommand{\rmat}{{\cal R}}
\newcommand{\oalg}{{\cal A}}
\newcommand{\falg}{{\cal F}}
\newcommand{\eich}{{\cal G}}
\newcommand{\hilb}{{\cal H}}
\newcommand{\calm}{{\cal M}}
\newcommand{\mod}{{\cal M}}

\newcommand{\bigrho}{\rho_\oplus}
\newcommand{\bigphi}{\phi_\oplus}

\newenvironment{bew}{Proof:}{\hfill$\Box$}
\newtheorem{bem}{Remark}
\newtheorem{bsp}{Example}
\newtheorem{axiom}{Axiom}
\newtheorem{de}{Definition}
\newtheorem{hypo}[de]{Hypothesis}
\newtheorem{satz}{Proposition}
\newtheorem{lemma}[satz]{Lemma}
\newtheorem{kor}[satz]{Corollary}
\newtheorem{theo}[satz]{Theorem}

\newcommand{\sbegin}[1]{\small\begin{#1}}
\newcommand{\send}[1]{\end{#1}\normalsize}

\sloppy

\newcommand{\lang}{1}

\title{The reduced Birman-Wenzl algebra of Coxeter type B}
\author{Reinhard H\"aring-Oldenburg\\ 
          Mathematisches Institut\\ Bunsenstr. 3-5\\
 37073 G\"ottingen, Germany\\
email: haering@cfgauss.uni-math.gwdg.de}
\date{17.7.97}
\maketitle

\abstract{We introduce a reduced form of a Birman-Murakami-Wenzl
Algebra associated to the braid group of Coxeter type B and investigate
its semisimplicity, Bratteli diagram and Markov trace. Applications in
knot theory and physics are outlined.
}

\section{Introduction}
To every Coxeter diagram a braid group is associated that has the same 
presentation as the Coxeter group but without the degree 2 relations
for the generators.  The braid group $\ZB_n$ of Coxeter type B has 
generators $\tau_i,i=0,1,\ldots n-1$.
Generators $\tau_i,i\geq1$  satisfy the relations of Artin's braid group 
(which is the braid group of Coxeter type A): 
\begin{eqnarray}
\tau_i\tau_j&=&\tau_j\tau_i \quad\mbox{if}\quad |i-j|>1\\
\tau_i\tau_j\tau_i&=&\tau_j\tau_i\tau_j\quad\mbox{if}\quad |i-j|=1
\end{eqnarray}
The generator $\tau_0$ has relations
\begin{eqnarray}
\tau_0\tau_1\tau_0\tau_1&=&\tau_1\tau_0\tau_1\tau_0\\
\tau_0\tau_i&=&\tau_i\tau_0\quad\mbox{if}\quad i\geq2
\end{eqnarray}
This braid group may be interpreted as the group of
symmetric braids or cylinder braids 
(see the graphical interpretation in section \ref{clsec}).

The group algebras of these braid groups  typically have lots 
of finite dimensional quotients.  The most important ones for 
Coxeter type A are Temperley-Lieb, Hecke and Birman-Murakami-Wenzl algebras.
Hecke algebras of arbitrary Coxeter type are already classics in this field.
Temperley-Lieb algebras of Coxeter type B have been introduced
by tom Dieck
in \cite{tD1} as algebras of symmetric tangles without crossings.

The standard Birman-Murakami-Wenzl algebra of type A imposes cubic relations
on its generators in a way that enables its interpretation  as an algebra
of tangles with a skein relation that comes from the Kauffman polynomial.

In full analogy a BMW algebra of Coxeter type B should be an extension by
an additional generator $Y$ related to $\tau_0$ which should satisfy
a cubic relation as well. It turns out, however, that such an algebra 
is rather intricate and deserves further study (see \cite{rho}).

In this	paper we define a reduced BMW algebra of type B where the additional
generator $Y$ satisfies a quadratic (Hecke type) relation. This may seem strange
at first but from the view of knot theory of B-type it is quite 
natural. Generalizations of this algebra where $Y$ may obey any
polynomial relation are considered in \cite{rhoak}.

We now outline the structure of the paper and point out the main results.
After a short review of the Birman-Wenzl algebra of A-type 
in section \ref{prelsec} we go on to define the reduced BMW algebra of B-type 
$\BB_n$
in section 	\ref{defsec} where a  number of fundamental relations
are established. They are used extensively in section \ref{wordsec} 
to determine  normal forms 	for words in $\BB_n$.
Un upper bound for the dimension is derived. 
Section \ref{connections} shows how to obtain 
the B-type Hecke algebra as a quotient of $\BB_n$.

Section \ref{clsec}  introduces the graphical interpretation 
of our algebra
and studies its classical limit. This will also give insight in the
relations chosen in the definition of $\BB_n$. 
The construction of a Markov trace fills section
\ref{trsec}. 

The main theorem of this paper is contained in section \ref{mainsec}.
We prove that $\BB_n$ is semisimple in the generic case and show how its 
simple components can be enumerated in terms of Young diagrams.
The Bratteli diagram is given and we show that the Markov trace is faithful.
 
T. tom Dieck has found a representation of $\BB_n$ on tensor product spaces.
In section \ref{repsec} we review his representation and show that it
allows to calculate the Markov trace as a matrix trace. 

The algebra $\BB_n$ has interesting applications both in physics and in knot
theory. They are outlined in the end of section \ref{repsec} and in section 
\ref{blinks}. 
The physical interest comes from the fact that the additional 
generator $Y$ may be interpreted as describing a boundary reflection
in a twodimensional quantum system.
The Markov trace allows to define an extension of the Kaufman polynomial
to links in the solid torus. 

A next goal would be to construct a tensor category \cite{rhobcat}
where $\BB_n$ is the endomorphism set of a $n$-fold tensor product
of a generating simple element.

{\em Acknowledgement}: Tammo tom Dieck deserves thanks for suggesting the 
study of this algebra and for many stimulating discussions. Thanks also
for the financial
support from the Deutsche Forschungsgemeinschaft.

\section{Preliminaries: The A-type BMW algebra}
\label{prelsec}

We review the definition of the Birman-Murakami-Wenzl algebra 
\cite{we1} in our notation	and collect a stock of relations 
that will be needed later on.

\begin{de}\label{defba}
Let $R$ denote an integral domain. Assume that $q,\lambda,x$ are
units in $R$ and define $\delta:=q-q^{-1}$. Assume that the relation 
\begin{equation} 
x\delta=\delta-\lambda+\lambda^{-1}
\end{equation}
holds. The Birman-Wenzl algebra of type A with $n$ strands $\BA_n(R)$ 
is defined as the algebra generated by invertible
$X_1,\ldots,X_{n-1}$.
The relations read:
\begin{eqnarray}
X_iX_j&=&X_jX_i\qquad |i-j|>1\label{def2}\\
X_iX_jX_i&=&X_jX_iX_j\qquad |i-j|=1\label{def3}\\
X_ie_i&=&e_iX_i=\lambda e_i\label{def4}\\
e_iX_{i-1}^{\pm1}e_i&=&\lambda^{\mp1}e_i\label{def5}
\\
e_i^2&=&xe_i\label{lem1a}\\
X_i^{-1}&=&X_i-\delta+\delta e_i\label{lem1d}\\
X_i^2&=&1+\delta X_i-\delta\lambda e_i\label{lem1e}\\
X_i^3&=&X_i^2(\lambda+\delta)+X_i(1-\lambda\delta)-\lambda\label{lem1b}\\
X_i^{-2}&=&1+\delta^2-\delta X_i+\delta(\lambda^{-1}-\delta)e_i
  = 1-\delta X_i^{-1}+\delta\lambda^{-1} e_i\label{lem1q}\\
0&=&(X_i-\lambda)(X_i+q^{-1})(X_i-q)\label{lem1c}\\
e_ie_j&=&e_je_i\qquad |i-j|>1\label{lem1f}\\
X_i^{-1}X_j^{\pm1}X_i&=&X_jX_i^{\pm1}X_j^{-1}\qquad |i-j|=1\label{lem1g}\\
e_iX_jX_i&=&X_j^\pm X_i^\pm e_j\qquad |i-j|=1\label{lem1h}\\
e_iX_j^{\pm1}e_i&=&\lambda^{\mp1}e_i\qquad |i-j|=1\label{lem1j}\\
e_ie_je_i&=&e_i\qquad |i-j|=1\label{lem1l}\\
X_i^{\pm1}e_je_i&=&X_j^{\mp1}e_i\qquad |i-j|=1\label{lem1k}\\
e_ie_jX_i^{\pm1}&=&e_iX_j^{\mp1}\qquad |i-j|=1\label{lem1kk}\\
e_iX^\pm_jX^\pm_i&=&e_ie_j\qquad |i-j|=1\label{lem1m}\\
X_i^\pm X_j^\pm e_i&=&e_je_i\qquad |i-j|=1\label{lem1mm}\\
X_ie_jX_i^{-1}&=&X_j^{-1}e_iX_j\qquad |i-j|=1\label{lem1i}\\
X_ie_jX_i&=&X_j^{-1}e_iX_j^{-1}\qquad |i-j|=1\label{lem1p}
\end{eqnarray}
\end{de}

\begin{lemma} \label{balemma}
If $\delta$ is invertible one may define  
\begin{equation}\label{edef}e_i:=1-\frac{X_i-X_i^{-1}}{\delta}
\end{equation}
and restrict the relations to (\ref{def2})-(\ref{def5}).
\end{lemma}
\begin{bew}
We have to show that the remaining relations are implied 
by this smaller set.
The proofs are mostly easy. We only comment on some of them. 
To show (\ref{lem1a}) one replaces one of the $e_i$
on the left hand side by its definition (\ref{edef}) and
applys (\ref{def4}). Relations (\ref{lem1d})-(\ref{lem1c}) are
succesive rewritings of (\ref{edef}).

(\ref{lem1g}): $X_iX_jX_i=X_jX_iX_j\Rightarrow
                X_jX_iX_j^{-1}=X_i^{-1}X_jX_i\Rightarrow
                X_iX_j^{-1}X_i^{-1}=X_j^{-1}X_i^{-1}X_j$

(\ref{lem1h}) follows from  (\ref{edef}) and (\ref{lem1g}).
To show (\ref{lem1l}) one replaces the $e_j$ in the middle
by its definition.

(\ref{lem1k}): $X_i^\pm e_je_i=X_j^\mp X_j^\pm X_i^\pm e_je_i
\stackrel{(\ref{lem1h})}{=}
X_j^\mp e_iX_j^\pm X_i^\pm e_i=\lambda^\pm X_j^\mp e_iX_j^\pm e_i
\stackrel{(\ref{lem1j})}{=}
X_j^\mp e_i$

(\ref{lem1m}): Using (\ref{lem1l}), (\ref{lem1k}) and 
         (\ref{lem1h}) we calculate
\begin{eqnarray*}
e_iX_j^\pm X_i^\pm &=&e_ie_je_iX^\pm _jX^\pm _i=
e_iX_i^{\mp}X_i^\pm e_je_iX^\pm _jX_i^\pm \\
&=&e_iX_i^{\mp}X_j^{\mp}e_iX^\pm _jX^\pm _i=
e_iX_i^{\mp}X_j^{\mp}X^\pm _jX^\pm _ie_j=e_ie_j
\end{eqnarray*}
\end{bew}

\section{The definition of the reduced B-type BMW algebra B}	\label{defsec}

In this section we define the reduced Birman-Murakami-Wenzl algebra
of Coxeter type B. The choice of the base ring needs special attention 
to avoid the algebra from being smaller than expected. 

\begin{de} Let $R$ be an integral domain of the kind described
in definition \ref{defba} with an additional unit
$q_0\in R$ and further elements $A,q_1\in R$.
The reduced Birman-Wenzl Algebra of Coxeter B type
with $n$ strands 
$\BB_n(R)$  is  generated by invertible
$Y$, $X_1,\ldots,X_{n-1}$
Using the notation from definition \ref{defba}
the relations are (\ref{def2}) to (\ref{def5}) and in addition:
\begin{eqnarray}
X_1YX_1Y&=&YX_1YX_1\label{def6}\\
Y^2&=&q_1Y+q_0\label{def7}\\
YX_1Ye_1&=&e_1\label{def8}\\
YX_i&=&X_iY\qquad i>1\label{def9}\\
e_1Ye_1&=&Ae_1\label{def10}
\end{eqnarray}
\end{de}

In the further development we assume that the algebra is
non degenerate in the sense that $e_1$ is non zero and
has a vanishing annulator ideal in $R$ and that
 $e_1,Ye_1$ are
linearly independent. Otherwise the algebra may not be simisimple.

\deu{Wir betrachten nun $Y$ betreffende Relationen und f\"uhren dazu folgende
Abk\"urzungen ein:}
\eng{We study now relations involving $Y$. The following shortcuts will
be useful:}
\begin{eqnarray}
Y'_i&:=&X_{i-1}X_{i-2}\cdots X_1YX_1\cdots X_{i-2}X_{i-1}\\
Y_i&:=&X_{i-1}X_{i-2}\cdots X_1YX_1^{-1}\cdots X_{i-2}^{-1}X_{i-1}^{-1}
\end{eqnarray}

\begin{lemma}\label{lemma2}
\begin{eqnarray}
Y^{-1}&=&q_0^{-1}Y-q_1q_0^{-1}\label{lem2a}\\
Y_i^2&=&q_1Y_i+q_0\label{lem2aa}\\
Y_i^{-1}&=&q_0^{-1}Y_i-q_1q_0^{-1}\label{lem2aaa}\\
0&=&[X_1YX_1Y,\{Y,e_1,X_1\}]\label{lem2b}\\
Y'_iY'_j&=&Y'_jY'_i\label{lem2f}\\
Y'_{i+1}X_i^{-1}&=&X_iY'_i\qquad Y_{i+1}X_i=X_iY_i\label{lem2g}\\
0&=&[Y_i,X_j]=[Y_i,e_j]\qquad j\neq i,i-1\label{lem2h}\\
0&=&[Y'_i,X_j]=[Y'_i,e_j]\qquad j\neq i,i-1\label{lem2hh}\\
e_i&=&e_iY_iX_iY_i=Y_iX_iY_ie_i\label{lem2c}\\
e_i&=&e_iY'_iX_iY'_i=Y'_iX_iY'_ie_i\label{lem2cs}\\
e_iY_ie_i&=&Ae_i\label{lem2ds}\\
X_iY_iX_iY_i&=&Y_iX_iY_iX_i\label{lem2zopf}\\
Y_ie_{i-1}&=&\lambda^{-1}q_0^{-1}Y_{i-1}e_{i-1}-q_1q_0^{-1}\lambda^{-1}e_{i-1}
\label{lem2q1}\\
e_{i-1}Y_i&=&\lambda(q_0^{-1}-\delta)e_{i-1}Y_{i-1}+
   \lambda(\delta A-q_1q_0^{-1})e_{i-1}\label{lem2q2}\\
X_iY_{i+1}&=&
Y_iX_i-\delta Y_i+\delta Y_ie_i+\delta Y_{i+1}\\
   &&+(\delta^2\lambda-\delta\lambda q_0^{-1})e_iY_i
            +(\delta\lambda q_1q_0^{-1}-\delta^2\lambda A)e_i\label{lem2q3}\\   
(1-q_0\delta)X_iY_ie_i&=&
e_i(q_1\lambda-q_0\delta\lambda A)+q_0Y_ie_i\label{lem2e}\\
e_{i-1}Y'_i&=&\lambda e_{i-1}Y_{i-1}'^{-1}\label{lem2m}\\
Y'_ie_{i-1}&=&\lambda Y_{i-1}'^{-1}e_{i-1}\label{lem2mm}\\
e_1Y_1X_2e_1&=&q_0e_1Y_1e_2e_1+q_1\lambda^{-1}e_1\label{lem2p}\\
Y_{i+1}Y_{i}&=&X_{i}Y_{i}X_{i}Y_{i}-\delta q_1X_{i}Y_{i}-\delta q_0X_{i}
 +\delta q_0^{-1} Y_{i}e_{i}Y_{i}-
 \delta q_1q_0^{-1}e_{i}Y_{i}\label{lem2z}
\end{eqnarray}
\end{lemma}
\begin{bew}
(\ref{lem2a}), (\ref{lem2aa}) \deu{und}\eng{and} (\ref{lem2aaa}) 
are verified easily.

(\ref{lem2b}): \deu{Mittels}\eng{Using} (\ref{def6}) 
\deu{ist zun\"achst}\eng{we have}
$X_1X_1YX_1Y=X_1YX_1YX_1$.
Hence  $X_1YX_1Y$ commutes with
$X_1$, and also with $X_1^{-1}$. 
But then, using (\ref{edef}), 
we see that it also commutes with $e_1$.

(\ref{lem2f}): $[Y,Y'_1]=[Y,Y'_2]=0$ \deu{ist trivial und f\"ur}
\eng{is trivial. For} $i>1$ \deu{gilt induktiv}\eng{the claim 
follows by induction}: 
$[Y,Y'_i]=0\Rightarrow[Y,Y'_{i+1}]=[Y,X_iY'_iX_i]=0$. 
\deu{F\"ur den allgemeinen Fall}\eng{In the general case} $[Y'_j,Y'_i]$ 
\deu{sei ObdA}\eng{we may assume}  $j<i$. 
\deu{Dann folgt der Induktionsschritt mit Hilfe von}
\eng{Then the induction step is shown using} (\ref{lem2h}): 
$[Y'_j,Y'_i]=[X_{j-1}Y'_{j-1}X_{j-1},Y'_i]=0$.

(\ref{lem2g}): \deu{Unmittelbare Konsequenzen der Definitionen.}\eng{trivial}

(\ref{lem2h},\ref{lem2hh}): \deu{F\"ur}\eng{For} $j\geq i+1$ 
\deu{folgt die Vertauschbarkeit aus}\eng{follows commutativity from} 
(\ref{def2},\ref{def9}) \deu{und f\"ur}\eng{and for} $j\leq i-1$ 
\deu{handelt es sich um eine einfache Anwendung der Relation}
\eng{it is an application of equation} (\ref{def3}).
\deu{Die Vertauschbarkeit mit $e_j$ folgt jeweils aus der mit}
\eng{Commutativity with $e_j$ follows from that with} $X_j$.

(\ref{lem2c},\ref{lem2cs}): \deu{Die Beweise erfolgen durch Induktion. 
Induktionsanfang ist}\eng{The proofs are by induction starting from} (\ref{def8}) 
\deu{und die gespiegelte Version}\eng{and its mirror version} 
$e_1=e_1YX_1Y$, \deu{die sich einfach beweisen l\"a\3t}
\eng{which may be proven easily}:
\[\lambda e_1YX_1Y=e_1X_1YX_1Y\stackrel{(\ref{lem2b})}{=}
X_1YX_1Ye_1=\lambda YX_1Ye_1=\lambda e_1\]

\deu{Induktionsschritt f\"ur}\eng{The induction step for } (\ref{lem2c}) 
\deu{benutzt}\eng{uses} (\ref{lem1h})\deu{, um}\eng{ to express} $e_{i+1}$
\deu{durch}\eng{in terms of} $e_i$\deu{ auszudr\"ucken}:
\begin{eqnarray*}
Y'_{i+1}X_{i+1}Y'_{i+1}e_{i+1}
&=&X_iY'_iX_iX_{i+1}X_iY'_iX_iX_i^{-1}X_{i+1}^{-1}e_iX_{i+1}X_i\\
&=&X_iY'_iX_{i+1}X_iX_{i+1}Y'_iX_{i+1}^{-1}e_iX_{i+1}X_i\\
&=&X_iX_{i+1}Y'_iX_iX_{i+1}X_{i+1}^{-1}Y'_ie_iX_{i+1}X_i\\
&=&X_iX_{i+1}Y'_iX_iY'_ie_iX_{i+1}X_i=X_iX_{i+1}e_iX_{i+1}X_i=e_{i+1}
\end{eqnarray*}
\deu{Induktionsschritt f\"ur}\eng{Induction step for} (\ref{lem2cs}):
\begin{eqnarray*}
e_{i+1}Y_{i+1}X_{i+1}Y_{i+1}
&=&e_{i+1}X_iY_iX_i^{-1}X_{i+1}X_iY_iX_i^{-1}
=e_{i+1}X_iY_iX_{i+1}X_iX_{i+1}^{-1}Y_iX_i^{-1}\\
&=&e_{i+1}X_iX_{i+1}Y_iX_iY_iX_{i+1}^{-1}X_i^{-1}\
=X_iX_{i+1}e_iY_iX_iY_iX_{i+1}^{-1}X_i^{-1}\\
&=&X_iX_{i+1}e_iX_{i+1}^{-1}X_i^{-1}=e_{i+1}\\
\end{eqnarray*}

(\ref{lem2ds}): 
\deu{Induktionsschritt}\eng{Induction step}:
\begin{eqnarray*}
e_iY_ie_i&=&e_iX_{i-1}Y_{i-1}X_{i-1}^{-1}e_i=
e_ie_{i-1}X_i^{-1}Y_{i-1}X_{i-1}^{-1}e_i\\
&=&e_ie_{i-1}Y_{i-1}X_i^{-1}X_{i-1}^{-1}e_i
=e_ie_{i-1}Y_{i-1}e_{i-1}X_i^{-1}X_{i-1}^{-1}\\
&=&Ae_ie_{i-1}X_i^{-1}X_{i-1}^{-1}\stackrel{(\ref{lem1kk})}{=}Ae_i
\end{eqnarray*}

(\ref{lem2zopf}): \deu{Der Induktionsanfang ist klar. Zum Induktionsschritt}
\eng{Again, the proof is by induction. The step is}:
\begin{eqnarray*}
Y_iX_iY_iX_i&=&X_{i-1}Y_{i-1}X_{i-1}^{-1}X_iX_{i-1}Y_{i-1}X_{i-1}^{-1}X_i\\
&=&X_{i-1}Y_{i-1}X_iX_{i-1}X_{i}^{-1}Y_{i-1}X_{i-1}^{-1}X_i\\
&=&X_{i-1}X_iY_{i-1}X_{i-1}Y_{i-1}X_i^{-1}X_{i-1}^{-1}X_i\\
&=&X_{i-1}X_iY_{i-1}X_{i-1}Y_{i-1}X_{i-1}X_i^{-1}X_{i-1}^{-1}\\
&=&X_{i-1}X_iX_{i-1}Y_{i-1}X_{i-1}Y_{i-1}X_i^{-1}X_{i-1}^{-1}\\
&=&X_{i}X_{i-1}X_{i}Y_{i-1}X_{i-1}X_i^{-1}Y_{i-1}X_{i-1}^{-1}\\
&=&X_{i}X_{i-1}Y_{i-1}X_{i}X_{i-1}X_i^{-1}Y_{i-1}X_{i-1}^{-1}\\
&=&X_{i}X_{i-1}Y_{i-1}X_{i-1}^{-1}X_{i}X_{i-1}Y_{i-1}X_{i-1}^{-1}\\
&=&X_{i}Y_{i}X_{i}Y_{i}
\end{eqnarray*}

(\ref{lem2q1}),(\ref{lem2q2}),(\ref{lem2q3}):
\begin{eqnarray*}
 Y_ie_{i-1}&=&X_{i-1}Y_{i-1}X_{i-1}^{-1}e_{i-1}=\lambda^{-1}X_{i-1}Y_{i-1}e_{i-1}\\
 &\stackrel{(\ref{lem2cs})}{=}&\lambda^{-1}Y_{i-1}^{-1}e_{i-1}
 =\lambda^{-1}q_0^{-1}Y_{i-1}e_{i-1}-q_1q_0^{-1}\lambda^{-1}e_{i-1}\\
e_{i-1}Y_i&=&e_{i-1}X_{i-1}Y_{i-1}X_{i-1}^{-1}=
 \lambda e_{i-1}Y_{i-1}X_{i-1}^{-1}\\
&=& \lambda e_{i-1}Y_{i-1}X_{i-1}-\delta \lambda e_{i-1}Y_{i-1}+
      \delta \lambda e_{i-1}Y_{i-1}e_{i-1}\\
&=&\lambda e_{i-1}Y_{i-1}^{-1}-\delta \lambda e_{i-1}Y_{i-1}+
      \delta \lambda A e_{i-1}\\
&=&q_0^{-1}\lambda e_{i-1}Y_{i-1}-q_1q_0^{-1}\lambda e_{i-1}-
   \delta \lambda e_{i-1}Y_{i-1}+\delta \lambda A e_{i-1}\\
&=&\lambda(q_0^{-1}-\delta)e_{i-1}Y_{i-1}+\lambda(\delta A-q_1q_0^{-1})e_{i-1}\\
X_iY_{i+1}&=&X_i^2Y_iX_i^{-1}=\\
&=&Y_iX_i^{-1}+\delta Y_{i+1}-\delta\lambda e_iY_iX_i^{-1}\\
&=&Y_iX_i-\delta Y_i+\delta Y_ie_i+\delta Y_{i+1}-\delta\lambda e_iY_iX_i
            +\delta^2\lambda e_iY_i-\delta^2\lambda e_iY_ie_i\\
&=&Y_iX_i-\delta Y_i+\delta Y_ie_i+\delta Y_{i+1}-\delta\lambda q_0^{-1}e_iY_i
            \\
			&&+\delta\lambda q_1q_0^{-1}e_i+\delta^2\lambda e_iY_i-
			   \delta^2\lambda A e_i\\
&=&Y_iX_i-\delta Y_i+\delta Y_ie_i+\delta Y_{i+1}+
                 (\delta^2\lambda-\delta\lambda q_0^{-1})e_iY_i
            \\&&+(\delta\lambda q_1q_0^{-1}-\delta^2\lambda A)e_i
\end{eqnarray*}

(\ref{lem2e}):
\begin{eqnarray*}
X_iY_ie_i&=&X_iY_iY_iX_iY_ie_i=q_1X_iY_iX_iY_ie_i+q_0 X_i^2Y_ie_i\\
&=&q_1X_ie_i+q_0(1+\delta X_i-\delta\lambda e_i)Y_ie_i\\
&=&q_1\lambda e_i+q_0Y_ie_i+q_0\delta X_iY_ie_i-q_0\delta\lambda Ae_i\\
\Rightarrow(1-q_0\delta)X_iY_ie_i&=&
e_i(q_1\lambda-q_0\delta\lambda A)+q_0Y_ie_i
\end{eqnarray*}

(\ref{lem2p}): \deu{Wir beweisen die folgende, \"aquivalente Relation}
\eng{We prove the following equivalent relation}:
\begin{eqnarray*}
e_1Y_1e_2e_1&=&e_1Y_1X_1X_2e_1=e_1Y_1X_1Y_1Y_1^{-1}X_2e_1=
e_1Y_1^{-1}X_2e_1\\
&=&q_0^{-1}e_1Y_1X_2e_1-q_1q_0^{-1}e_1X_2e_1
=q_0^{-1}e_1Y_1X_2e_1-q_1q_0^{-1}\lambda^{-1}e_1
\end{eqnarray*}

(\ref{lem2m},\ref{lem2mm}) \deu{beweist man nach dem Schema}
\eng{is proven according to the scheme }
\[e_{i-1}Y'_i=e_{i-1}X_{i-1}Y'_{i-1}X_{i-1}=
\lambda e_{i-1}Y'_{i-1}X_{i-1}Y'_{i-1}Y_{i-1}'^{-1}=
\lambda e_{i-1}Y_{i-1}'^{-1}\]

(\ref{lem2z}):
\begin{eqnarray*}
Y_{i+1}Y_i&=&X_iY_iX_i^{-1}Y_i\\
&=&X_iY_iX_iY_i-\delta X_iY_i^2+\delta X_iY_ie_iY_i\\
&=&X_iY_iX_iY_i-\delta q_1 X_iY_i-\delta q_0 X_i+\delta Y_i^{-1}e_iY_i\\
&=&X_iY_iX_iY_i-\delta q_1 X_iY_i-\delta q_0 X_i+\delta q_0^{-1} Y_ie_iY_i-
  \delta q_1q_0^{-1}e_iY_i
\end{eqnarray*}

\end{bew}

Our non degeneracy assumptions introduce relations
among the parameters. 

\begin{lemma} The assumption that $e_1$ has non vanishing annulator
ideal leads to the requirement
\begin{equation}\label{Avalue}
A(1-q_0\lambda)=q_1x
\end{equation}
The additional assumption that $Ye_1$ and $e_1$  are linearly independent
leads to the equation 
\begin{equation}\label{q0value}
q_0-q_0^{-1}=-\delta
\end{equation}
\end{lemma}
\begin{bew}
\begin{eqnarray*}
e_1Ye_1&=&e_1YYX_1Ye_1=q_1e_1YX_1Ye_1+q_0e_1X_1Ye_1
=q_1xe_1+q_0\lambda e_1Ye_1\\
\Rightarrow && (1-q_0\lambda)e_1Ye_1=q_1xe_1\quad
\Rightarrow\quad A(1-q_0\lambda)=q_1x
\end{eqnarray*}
To obtain the second relation we observe that (\ref{def8})
implies $Ye_1=X_1^{-1}Y^{-1}e_1$. We multiply by $q_0$ and calculate
\begin{eqnarray*}
q_0Ye_1&=&q_0X_1^{-1}Y^{-1}e_1=(X_1-\delta+\delta e_1)(Y-q_1)e_1\\
&=&  X_1Ye_1  -q_1X_1e_1  -\delta Ye_1  +\delta q_1e_1
+\delta e_1Ye_1-\delta q_1e_1^2\\
&=&(q_0^{-1}Y-q_0^{-1}q_1)e_1-q_1\lambda e_1-\delta Ye_1+\delta q_1e_1
+\delta A e_1-\delta q_1xe_1\\
&=&e_1(-q_1q_0^{-1}-q_1\lambda+\delta q_1+\delta A-\delta q_1 x)+
Ye_1(q_0^{-1}-\delta)
\end{eqnarray*}
The coefficient of $Ye_1$ is (\ref{q0value}). The coefficient of $e_1$
vanishes when (\ref{Avalue}) and (\ref{q0value}) hold.
\end{bew}

From now on we will always assume that these relations hold
in the ground ring. 

Using the relations of lemma \ref{lemma2} one sees that the ideal
generated by $e_1$ in $\BB_2$ is spanned by $e_1,Ye_1,e_1Y,Ye_1Y$. 
Using the relations of the above lemma one may 
(by construction of a twodimensional irreducible representation) show 
that the ideal is indeed four dimensional and
hence that the nondegeneracy assumptions imply 
no further relations among the parameters.
We don't go into details of this but see \cite{rhoak}
for a detailed exposition of such arguments in a more
complicated case.

At this stage of the development it is useful to look ahead
to the classical limit of the algebra we shall discuss later on. 
Such a limit should have $X_1=X_1^{-1}$ which is implied by
$q\rightarrow1$. Furthermore, one would expect that $Y$
as well should obey a Coxeter relation $Y^2=1$ in the limit. 
It is therefore
reasonable to choose \begin{equation}q_0=q^{-1}\end{equation}
among the solutions of (\ref{q0value}) as we will do from now on.

The generic ground ring that we will use is:
\begin{de}
The ring $R_0$ is defined to be the quotient of the polynomial ring  
$\CC[q,q^{-1},q_0,q_0^{-1}, \delta,\delta^{-1},
\lambda,\lambda^{-1},q_1,A]$ quotiented by the relations 
(\ref{Avalue}), $\delta=q-q^{-1}$ and the Laurent style relations
$qq^{-1}=1$ and so on. Its quotient field is denoted by $K_0$.
\end{de}
Here we have already eliminated $q_0$.
In the quotient ring of $R_0$ we can solve the
equations defining the ideal uniquely. Hence this ideal
is primary and therfore $R_0$  is an integral domain.
Therefore $R_0$ is embedded in $K_0$.

\begin{bem} \label{involution} 
The algebra $\BB_n$ has an involution given by
\begin{equation}
X_i^\ast:=X_i^{-1},Y^\ast:=Y^{-1},
q^\ast:=q^{-1},\lambda^\ast:=\lambda^{-1},
q_0^\ast:=q_0^{-1}, q_1^\ast:=-q_1q_0^{-1}\end{equation}. 
This implies $\delta^\ast=-\delta,e_i^\ast=e_i,
A^\ast:=(A-q_1x)/q_0$.

A second involution $a\mapsto\overline{a}$ exists that
fixes all parameters and generators.
\end{bem}

\section{The word problem in $\BB_n$}	\label{wordsec}

In this section we single out a set of words in standard form
that linearly generate $\BB_n$. Although this 
does not lead to a linear basis of $\BB_n$, it allows to determine
a tight upper bound for the dimension.

\begin{satz}\label{wortsatz}
\eng{Every element in $\BB_n$ is a linear combination of words of the form}
$w_1\gamma w_2$, \deu{wobei}\eng{where} $w_i\in\BB_{n-1}$ \deu{und}\eng{and}
$\gamma\in\Gamma_n:=\{1,e_{n-1},X_{n-1},Y_{n}\}$
\end{satz}
\begin{bew} \deu{Induktion nach n.}\eng{We prove the proposition by induction.}
\eng{The case $n=1$ is trivial and $n=2$ can also be verified easily.}

\eng{Let}\deu{Sei} $w_0\gamma_0w_1\gamma_1\cdots w_k\gamma_kw_{k+1}\in\BB_n$ 
\deu{ein beliebiges Wort}\eng{be an arbitrary word}.
\deu{Es reicht zu zeigen, da\3 je zwei benachbarte}
\eng{It suffices to show that any two neighbouring} $\gamma_i$ 
\deu{zusammengefa\3t werden k\"onnen.}
\eng{can be combined together.}
\deu{Die relevante Situation ist also}
\eng{Hence the situation we have to investigate is } $w=\gamma_1w_1\gamma_2,
w_1\in\BB_{n-1},\gamma_1,\gamma_2\in\Gamma_n$. 
\deu{Nach Induktionsannahme ist}\eng{By inducton hypothesis we have}
$w_1=u_1\alpha u_2,u_i\in\BB_{n-2},\alpha\in\Gamma_{n-1}$ 
\deu{und daher}\eng{and hence}
$w=\gamma_1u_1\alpha u_2\gamma_2=u_1\gamma_1\alpha\gamma_2 u_2$. 
\deu{Es reicht also,}\eng{Thus it suffices to investigate}
$w'=\gamma_1\alpha\gamma_2$\deu{ zu betrachten}. 
\deu{Weiter sind die F\"alle}\eng{The cases}
$\gamma_1=1$ \deu{oder}\eng{or} $\gamma_2=1$ \eng{are} trivial. 
\deu{Wir betrachten jetzt sukzessive die vier m\"oglichen Werte von $\alpha$.}
\eng{We now investiagte in turn the four possible values of $\alpha$.}

1. \deu{Fall}\eng{Case} $\alpha=1$: 
\deu{Die folgende Tabelle notiert die Relation in der Algebra,
die es gestattet, das Produkt $\gamma_1\gamma_2$ auf die im Satz angegebene
Form zu bringen.}
\eng{The following table gives the relation that allows to reduce the product
$\gamma_1\gamma_2$ to the standard form of the proposition.}
\[\begin{array}{c||c|c|c}
\gamma_1\backslash\gamma_2 & Y_n & e_{n-1} & X_{n-1}\\\hline\hline
Y_n &(\ref{lem2aa}) &(\ref{lem2q1}) & (\ref{lem2g})\\\hline
e_{n-1} &(\ref{lem2q2}) & (\ref{lem1a})&(\ref{def4})\\\hline
X_{n-1} & (\ref{lem2q3})& (\ref{def4}) & (\ref{lem1e})
\end{array}\]

2. \deu{Fall}\eng{Case} $\alpha=X_{n-2}$: 
\[\begin{array}{c||c|c|c}
\gamma_1\backslash\gamma_2 & Y_n & e_{n-1} & X_{n-1}\\\hline\hline
Y_n & =X_{n-2}Y_n^2 \quad (\ref{lem2aa})&=X_{n-2}Y_ne_{n-1}\quad(\ref{lem2q1})
   &=X_{n-2}Y_nX_{n-1}\quad (\ref{lem2g})\\\hline
e_{n-1} &=e_{n-1}Y_nX_{n-2}\quad(\ref{lem2q2}) &
  (\ref{def5}) &(\ref{lem1m})\\\hline
X_{n-1} & (\ref{lem2q3})& (\ref{lem1mm}) & (\ref{def3})
\end{array}\]

3. \deu{Fall}\eng{Case} $\alpha=e_{n-2}$: 
\[\begin{array}{c||c|c|c}
\gamma_1\backslash\gamma_2 & Y_n & e_{n-1} & X_{n-1}\\\hline\hline
Y_n & =e_{n-2}Y_n^2\quad (\ref{lem2aa})&(\ref{lem2q1}) & (\ref{lem2g})\\\hline
e_{n-1} &(\ref{lem2q2}) & (\ref{lem1l})&(\ref{lem1kk})\\\hline
X_{n-1} & (\ref{lem2q2})& (\ref{lem1k}) & (\ref{lem1p})
\end{array}\]

4. \deu{Fall}\eng{Case} $\alpha=Y_{n-1}$:
\deu{Hier sind zum Teil etwas aufwendigere Rechnungen erforderlich, die im Anschlu\3
an die Tabelle wiedergegeben werden.}
\eng{This case requiers more complex calculations which are given below.} 
\[\begin{array}{c||c|c|c}
\gamma_1\backslash\gamma_2 & Y_n & e_{n-1} & X_{n-1}\\\hline\hline
Y_n & (\ref{flaa}) &(\ref{flab}) & (\ref{flac})\\\hline
e_{n-1} &(\ref{flba}) & (\ref{lem2ds})&(\ref{flbc})\\\hline
X_{n-1} &(\ref{flca}) & \mbox{\deu{analog}\eng{analog.} }(\ref{flbc})& (\ref{flcc}) 
\end{array}\]

\begin{eqnarray}
Y_nY_{n-1}Y_n&=&X_{n-1}Y_{n-1}X_{n-1}^{-1}Y_{n-1}X_{n-1}Y_{n-1}X_{n-1}^{-1}\label{flaa}\\
&\stackrel{(\ref{lem2zopf})}{=}&
  X_{n-1}Y_{n-1}Y_{n-1}X_{n-1}Y_{n-1}X_{n-1}^{-1}X_{n-1}^{-1}\nonumber\\
&=&q_1 X_{n-1}Y_{n-1}X_{n-1}Y_{n-1}X_{n-1}^{-2}+
   q_0 X_{n-1}^2Y_{n-1}X_{n-1}^{-2}\nonumber\\
&=&q_1 Y_{n-1}X_{n-1}Y_{n-1}X_{n-1}^{-1}+
   q_0 X_{n-1}^2Y_{n-1}X_{n-1}^{-2}\nonumber\\
&=&q_1 Y_{n-1}Y_n\nonumber\\
  &&+ q_0 (1+\delta X_{n-1}-\delta\lambda e_{n-1})Y_{n-1}
     (1+\delta^2-\delta X_{n-1}+\delta(\lambda^{-1}-\delta)e_{n-1})\nonumber\\
&&    \mbox{\deu{Damit ist der Ausdruck auf die anderen F\"alle der Tabelle
 reduziert.}\eng{This reduces the problem to the other cases.}}\nonumber
\end{eqnarray}

\begin{eqnarray}\label{flab}
\lefteqn{Y_nY_{n-1}e_{n-1}\stackrel{(\ref{lem2z})}{=}}\\
&=&X_{n-1}Y_{n-1}X_{n-1}Y_{n-1}e_{n-1}-\delta q_1 X_{n-1}Y_{n-1}e_{n-1}-
\delta q_0 X_{n-1}e_{n-1}\nonumber\\
&&+\delta q_0^{-1} Y_{n-1}e_{n-1}Y_{n-1}e_{n-1}-
\delta q_1q_0^{-1}e_{n-1}Y_{n-1}e_{n-1}\nonumber\\
&=&X_{n-1}e_{n-1}-\delta q_1 Y_{n-1}^{-1}e_{n-1}-\delta q_0 \lambda e_{n-1}\nonumber\\
&&+\delta q_0^{-1}A Y_{n-1}e_{n-1}-\delta q_1q_0^{-1}Ae_{n-1}\nonumber\\
&=&\lambda e_{n-1}-\delta q_1 Y_{n-1}^{-1}e_{n-1}-
\delta q_0 \lambda e_{n-1}
+\delta q_0^{-1}A Y_{n-1}e_{n-1}-\delta q_1q_0^{-1}Ae_{n-1}\nonumber
\end{eqnarray}

\begin{eqnarray}\label{flac}
Y_nY_{n-1}X_{n-1}&\stackrel{(\ref{lem2z})}{=}&
X_{n-1}Y_{n-1}X_{n-1}Y_{n-1}X_{n-1}-\delta q_1 X_{n-1}Y_{n-1}X_{n-1}\\
&&-\delta q_0 X_{n-1}X_{n-1}
+\delta q_0^{-1} Y_{n-1}e_{n-1}Y_{n-1}X_{n-1}\nonumber\\
&&-\delta q_1q_0^{-1}e_{n-1}Y_{n-1}X_{n-1}\nonumber\\
&=&
Y_{n-1}X_{n-1}Y_{n-1}X_{n-1}^2-\delta q_1 Y_nX_{n-1}^2-
\delta q_0 X_{n-1}^2\nonumber\\
&&+\delta q_0^{-1} Y_{n-1}e_{n-1}Y_{n-1}^{-1}-
\delta q_1q_0^{-1}e_{n-1}Y_{n-1}^{-1}\nonumber\\
&& \mbox{\deu{Nur die ersten beiden Terme sind noch nicht reduziert.}
\eng{Only the first and second term are not yet reduced.}}\nonumber\\
Y_{n-1}X_{n-1}Y_{n-1}X_{n-1}^2&=&Y_{n-1}Y_nX_{n-1}^3\nonumber\\
&&\mbox{\deu{Das reduziert man mittels}
\eng{This is reduced using}(\ref{lem1b},\ref{lem2q1})}\nonumber\\
Y_nX_{n-1}^{-2}&=&Y_n(1+\delta^2)-\delta Y_nX_{n-1}+
\delta(\lambda^{-1}-\delta)Y_ne_{n-1}\nonumber\\
&=&Y_n(1+\delta^2)-\delta X_{n-1}Y_{n-1}+
\delta(\lambda^{-1}-\delta)Y_ne_{n-1}\nonumber\\
&&\mbox{\deu{Das reduziert sich mit}\eng{This can be reduced using} 
(\ref{lem2q1})}\nonumber
\end{eqnarray}

\begin{equation}
\label{flba}
e_{n-1}Y_{n-1}Y_n=e_{n-1}Y_{n-1}X_{n-1}Y_{n-1}X_{n-1}^{-1}=\lambda^{-1}e_{n-1}
\end{equation}

\begin{equation}\label{flbc} e_{n-1}Y_{n-1}X_{n-1}=
e_{n-1}Y_{n-1}X_{n-1}Y_{n-1}Y_{n-1}^{-1}=
e_{n-1}Y_{n-1}^{-1}
\end{equation}

\begin{equation}\label{flca} X_{n-1}Y_{n-1}Y_n=X_{n-1}Y_{n-1}X_{n-1}Y_{n-1}X_{n-1}^{-1}
=Y_{n-1}X_{n-1}Y_{n-1}
\end{equation}

\begin{eqnarray}\label{flcc}
X_{n-1}Y_{n-1}X_{n-1}&=&Y_nX_{n-1}^2=Y_n+\delta Y_nX_{n-1}-\delta\lambda Y_ne_{n-1}\\
&=&Y_n+\delta X_{n-1}Y_{n-1}-\delta\lambda Y_ne_{n-1}\nonumber\\
&&\mbox{\deu{Der letzte Term reduziert sich mit}
\eng{The last term can be  reduced using} (\ref{lem2q1})} \nonumber
\end{eqnarray}
\end{bew}

\eng{This shows that $\BB_n$ is finite dimensional.}

\begin{bem}\label{wortinv} \deu{Es ist offensichtlich, da\3 analoge
S\"atze gelten, wenn in} 
\eng{It is obvious that similar propositions hold if} 
\deu{$Gamma_n$} $Y_n$ \deu{oder}\eng{or} $X_{n-1}$ 
\deu{oder beide}\eng{or both in $\Gamma_n$}
\deu{durch ihr Inverses ersetzt werden}\eng{are replaced by their inverses}.
\end{bem}

\deu{Der folgende Satz ergibt sich jetzt als Korollar: }
\begin{satz} \label{yswort}
\deu{Die Aussage von Satz}\eng{In proposition} \ref{wortsatz} 
\deu{gilt auch wenn}\eng{one may replace} $\Gamma_n$ 
\eng{by}\deu{ersetzt wird durch} $\Gamma'_n:=\{1,e_{n-1},X_{n-1},Y'_n\}$.
\end{satz}
\begin{bew}
\deu{Es reicht zu zeigen, da\3}\eng{It suffices to show that} $Y_n$ 
\deu{durch Normalformworte mit}\eng{can be expressed using words in
normal form with } $Y'_n$\deu{ ausgedr\"uckt werden kann}. 
\deu{F\"ur}\eng{For} $n=1$ \deu{ist das trivial}\eng{this is trivial}. 
\deu{Induktion: In}\eng{Induction step: Express $Y_n$ in}
$Y_{n+1}=X_nY_nX_n^{-1}$ \deu{werde $Y_n$ durch Normalformworte ausgedr\"uckt.}
\eng{in terms of normal form words.}
\deu{Solange diese mit}\eng{If they are build with}  $1,X_{n-1}$ 
\deu{oder}\eng{or} $e_{n-1}$ \deu{als}\eng{as} $\gamma$ 
\deu{gebildet werden, ist nichts zubeweisen.}
\eng{there is nothing to show.}
\deu{Es bleibt also nur}\eng{The only  remaining case is:}
\begin{eqnarray*}
X_nY'_nX_n^{-1}&=&X_nY'_nX_nX_n^{-2}=
          Y'_{n+1}(1-\delta X_n^{-1}+\delta\lambda^{-1} e_n)\\
          &=& Y'_{n+1}-\delta Y'_{n+1}X_n^{-1}+
          \delta\lambda^{-1}Y'_{n+1}e_n 
           \stackrel{(\ref{lem2mm})}{=} Y'_{n+1}-\delta X_nY'_{n}+
          \delta Y_{n}'^{-1}e_n
\end{eqnarray*}
This shows that terms of this kind can be brought 
to the normal form as well.
\end{bew}

The  aim of the rest of this section is to determine an upper bound
for the dimension of $\BB_n$.

\begin{lemma}
$\BB_n$ is spanned linearly by the set $S_n$ defined recursively by:
\begin{eqnarray*}
S_1&:=&\{1,Y\}\\
S_n&:=&\Gamma'_1\cdots\Gamma'_nS_{n-1}
\end{eqnarray*}
More strongly, of the elements of $\Gamma'_1\cdots\Gamma'_n$ only
those of the following form are needed.
\[ Y'_iX_i\cdots X_je_{j+1}\cdots e_n,\qquad
X_i\cdots X_je_{j+1}\cdots e_n\]
Here $1\leq i\leq n$ and $i-1\leq j\leq n$. Thus the strings of $X$ and $e$
may be empty.
\end{lemma}
\begin{bew}
Proposition \ref{yswort} yields the following decomposition of
$\BB_n$ which implies the claim:
\begin{eqnarray*}
\BB_n&=&\BB_{n-1}\Gamma'_n\BB_{n-1}\\
&=&\BB_{n-2}\Gamma'_{n-1}\BB_{n-2}\Gamma'_n\BB_{n-1}=
\BB_{n-2}\Gamma'_{n-1}\Gamma'_n\BB_{n-1}\\
&=&\Gamma'_1\cdots\Gamma'_n\BB_{n-1}
\end{eqnarray*}
To show the second statement assume that $Y'_j$ appears
in the middle of a chain $Z_i\cdots Z_{j-1}Y'_jZ_{j+1}\cdots Z_n$
where $Z_s\in\Gamma'_s$. Then $Z_i\cdots Z_{j-1}$ commutes
with the rest of the chain and thus can be absorbed in the
right $\BB_{n-1}$.
Similarly, assume that there appears a $e_iX_{i+1}$ in such a chain.
Then one can rewrite this as $e_iX_{i+1}=e_ie_{i+1}X_i^{-1}$ and
now the $X_i^{-1}$ can be absorbed in the right $\BB_{n-1}$.
Thus all $X$ must appear to the left of all $e$ in the chain. 
This completes the proof of the given form.
\end{bew}

\begin{satz}\label{dimsatz}
There is a basis of $\BB_n$ consisting of elements of the form
$\alpha\beta\gamma$ where $\alpha$ is a product of $Y'$, 
$\gamma$ is a product of ${Y'}^{-1}$ and 
$\beta$ is an element of a basis of the A-type algebra $\BA_n$.
Together $\alpha$ and $\gamma$ contain at most $n$ factors 
$Y',Y'^{-1}$.

The dimension of $\BB_n$ is $\leq 2^n(2n-1)!!$.
\end{satz}
\begin{bew}
The proof is by induction on $n$. For $n=1$ it is trivial.
Now, assume the claim is already shown for $n-1$. 
To show the first statement it suffices to show
that we can move all $Y'_i$ that appear on the left hand side of 
our basis of $\BB_{n-1}$ through the outer $\Gamma'$ chain
to the left or, alternatively,  even to the right of $\BB_{n-1}$.
We investigate the various arising cases.
First assume that we have $e_{n-1}Y'_{n-1}$. Then we rewrite this as 
\begin{eqnarray*}
e_{n-1}Y'_{n-1}&=&
e_{n-1}Y'_{n-1}X_{n-1}Y'_{n-1}{Y'_{n-1}}^{-1}X_{n-1}^{-1}
=e_{n-1}{Y'_{n-1}}^{-1}X_{n-1}^{-1}\\
&=&\lambda e_{n-1}X_{n-1}^{-1}{Y'_{n-1}}^{-1}X_{n-1}^{-1}
=\lambda e_{n-1}{Y'_n}^{-1}
\end{eqnarray*}
If we have $e_ie_{i+1}Y'_i=e_iY'_ie_{i+1}$	we may 
apply the same reasoning twice to obtain $Y'_{i+2}e_ie_{i+1}$.
The remainig cases are such that we have 
$X_iY'_i=Y'_{i+1}X_i^{-1}=Y'_{i+1}X_i-\delta Y'_{i+1}+
\delta Y'_{i+1}e_i$. The first summand is of the desired
form.  In the second there may be a chain of $X$ left to the
 $Y'_{i+1}$ which may be commuted to the right and absorbed in the
$\BB_{n-1}$. The third summand is either of the desired form, 
or it may violate the
rule that no $e_i$ should appear in a chain on the left of a $X_i$.
But if this rule is violated, it may be restored by the same
argument as in the proof of the previous lemma.

None of our rewritings did change the number of $Y'$ and so we can't
have more than $n$ of them, at most one coming from each recursion
in the construction of $S_n$. By induction assumption the dimension 
of $\BB_{n-1}$ is less than $2^{n-1}(2n-3)!!$ and we have brought the 
$Y'$ safely outside the region of $\BA_n$ elements. 
From the theory of $\BA_n$
it follows that $2n-1$ different chains  $Z_i\cdots Z_n,
Z_j\in\{e_{i-1},X_{i-1}\}$ are needed. Each of these chains may
have a $Y'_i$ at its front. Hence we conclude that the dimensions increases 
at most by a factor $2(2n-1)$. Thus the claim follows.
\end{bew}

\section{Relation to the B type Hecke algebras}\label{connections}

\begin{de}							\label{heckedef}
Let $\HB_n$ denote the Hecke algebra of Coxeter type B 
\deu{Sie hat Generatoren}\eng{with generators}
$X_0,X_1,\ldots,X_{n-1}$ \deu{und Parameter}\eng{and parameters} $Q,Q_0$ 
\deu{sowie die Relationen}\eng{and relations}:
\begin{eqnarray}
X_0X_1X_0X_1&=&X_1X_0X_1X_0\label{hdef1}\\
X_iX_j&=&X_jX_i\qquad |i-j|>1\label{hdef2}\\
X_iX_jX_i&=&X_jX_iX_j\qquad |i-j|=1\label{hdef3}\\
X_i^2&=& (Q-1)X_i+Q\quad i\geq 0\label{hdef4}\\
X_0^2&=& (Q_0-1)X_0+Q_0
\end{eqnarray}
\end{de}

\begin{lemma}					\label{heckelemma}
\eng{Let} $I_n$ \deu{sei das von}\eng{be the ideal generated by} $e_{n-1}$ 
\deu{erzeugte Ideal} in $\BB_n$.
\deu{Jedes beliebige}\eng{Every other} $e_i$ 
\deu{erzeugt dasselbe Ideal und der Quotient ist isomorph zur Hecke-Algebra}
\eng{generates the same ideal and the quotient algebra is isomorphic to} 
$HB_n$.
\end{lemma}
\begin{bew}
The first relation follows from (\ref{lem1i})
\eng{which allows to express any $e_i$ in terms of any other $e_j$.}
The isomorphism $\BB_n/I_n\rightarrow\HB_n$ is given by
$X_i\mapsto q^{-1}X_i, Q=q^2, Y\mapsto -X_0q^{-1}(qq_1+\sqrt{4q+q^2q_1^2})/2,
2Q_0=2+qq_1^2-q_1\sqrt{4q+q^2q_1^2}$
\end{bew}

Of course one can avoid square roots by using a different normalization 
of the generators.

\begin{lemma} \label{ilem} $I_n=\BB_{n-1}e_{n-1}\BB_{n-1}$
\end{lemma}
\begin{bew}
The ideal is defined to be $I_n=\BB_{n}e_{n-1}\BB_{n}$.
If we apply proposition \ref{wortsatz} we obtain
\begin{eqnarray*}
I_n&=&\BB_{n-1}\Gamma'_n\BB_{n-1}e_{n-1}\BB_{n-1}\Gamma'_n\BB_{n-1}\\
&=&\BB_{n-1}\Gamma'_n\BB_{n-2}\Gamma'_{n-1}\BB_{n-2}
e_{n-1}\BB_{n-2}\Gamma'_{n-1}\BB_{n-2}\Gamma'_n\BB_{n-1}\\
&=&\BB_{n-1}\Gamma'_n\Gamma'_{n-1}
e_{n-1}\BB_{n-2}\Gamma'_{n-1}\Gamma'_n\BB_{n-1}
\end{eqnarray*}
Hence it suffices to establish that 
$\Gamma'_n\Gamma'_{n-1}e_{n-1}\subset\BB_{n-1}e_{n-1}$. 
This is done easily using the relations from lemma
\ref{balemma} and \ref{lemma2}.
\end{bew}

\section{Graphical Interpretation and the classical limit}
\label{clsec}

The definition of $\BB_n$ is inspired by  
B type knot theory. This section supplies the precise definition
of the graphical version of the algebra.

\begin{figure}[ht]
\begin{picture}(130,60)
\put(59,50){\mbox{$Y^{(G)}\equiv X_0^{(G)}$}}

\linethickness{0.2mm}
\put(2,50){\mbox{$\cdots$}}
\put(9,56){\mbox{{\small -3}}}
\put(14,56){\mbox{{\small -2}}}
\put(19,56){\mbox{{\small -1}}}
\put(40,56){\mbox{{\small 3}}}
\put(35,56){\mbox{{\small 2}}}
\put(30,56){\mbox{{\small 1}}}

\put(10,45){\line(0,1){10}}
\put(15,45){\line(0,1){10}}
\put(20,55){\line(1,-1){10}}
\put(20,45){\line(1,1){4}}
\put(26,51){\line(1,1){4}}
\put(35,45){\line(0,1){10}}
\put(40,45){\line(0,1){10}}
\put(42,50){\mbox{$\cdots$}}

\linethickness{0.4mm}
\put(90,45){\line(0,1){3}}
\put(90,50){\line(0,1){5}}
\linethickness{0.2mm}
\put(88,51){\oval(4,4)[l]}
\put(88,49){\line(1,0){5}}
\put(93,47){\oval(4,4)[tr]}
\put(93,55){\oval(4,4)[br]}
\put(95,56){\mbox{{\small 1}}}
\put(90,56){\mbox{{\small 0}}}

\put(100,45){\line(0,1){10}}

\put(107,50){\mbox{$\cdots$}}


\put(72,27){\mbox{$X_i^{(G)}$}}

\put(0,32){\line(1,-1){10}}
\put(0,22){\line(1,1){4}}
\put(6,28){\line(1,1){4}}
\put(45,32){\line(1,-1){10}}
\put(45,22){\line(1,1){4}}
\put(51,28){\line(1,1){4}}
\put(20,22){\line(0,1){10}}
\put(30,22){\line(0,1){10}}
\put(13,27){\mbox{$\cdots$}}
\put(35,27){\mbox{$\cdots$}}
\linethickness{0.4mm}
\put(90,22){\line(0,1){10}}
\linethickness{0.2mm}
\put(95,22){\line(0,1){10}}
\put(98,27){\mbox{$\cdots$}}
\put(102,32){\line(1,-1){10}}
\put(102,22){\line(1,1){4}}
\put(108,28){\line(1,1){4}}


\put(72,5){\mbox{$e_i^{(G)}$}}
\linethickness{0.4mm}
\put(90,0){\line(0,1){10}}
\linethickness{0.2mm}
\put(95,0){\line(0,1){10}}
\put(98,5){\mbox{$\cdots$}}
\put(106,0){\oval(8,8)[t]}
\put(106,10){\oval(8,8)[b]}
\put(20,0){\line(0,1){10}}
\put(30,0){\line(0,1){10}}
\put(12,5){\mbox{$\cdots$}}
\put(32,5){\mbox{$\cdots$}}
\put(5,0){\oval(8,8)[t]}
\put(5,10){\oval(8,8)[b]}
\put(45,0){\oval(8,8)[t]}
\put(45,10){\oval(8,8)[b]}

\end{picture}
\caption{\label{generat} 
The graphical interpretation of the generators as symmetric tangles (on the left)
and as cylider tangles (on the right)}

\end{figure}
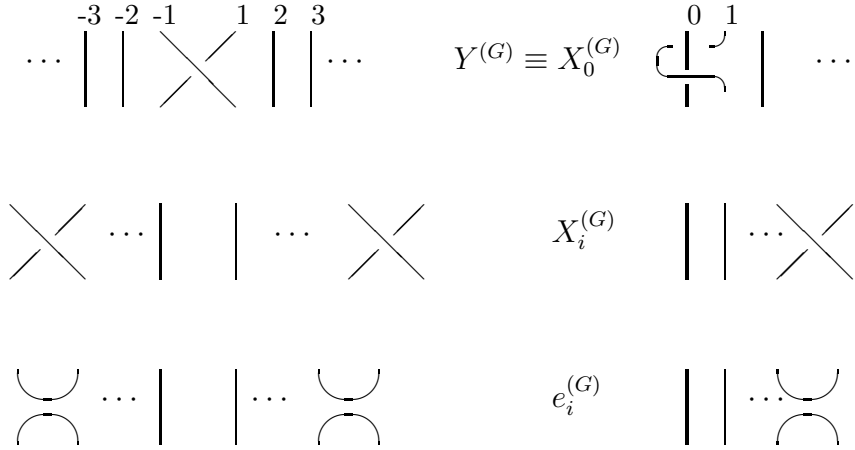

Let $R$ be an integral domain. Consider the free $R$ algebra generated by 
isotopy classes of ribbons in  $(\RR^2-\{0\})\times [0,1]$ 
between $n$ upper and $n$ lower intervals imbedded on
the line $\RR^+\times 0\times 1$ resp. $\RR^+\times 0\times 1$.
There may be ribbon components that are not connected to 
these endpoints. Multiplication is given by putting 
the graphs on top of each other.
Next, restrict the attention to the subalgebra that consists of
thoses isotopy calsses that have a representation as a product of the
generators $X_i^{(G)},e_i^{(G)},Y^{(G)},1\leq i\leq n-1$ 
from figure \ref{generat}. 
We define $G\BB_n(R)$ 
(where $R$ is  as in the definition of $\BB_n$ with (for the moment) $\delta$ 
invertible)
to be the quotient of this algebra by the relations
(\ref{def4}), (\ref{def5}), (\ref{def7}), (\ref{def10}). 
The remaining relations	in the definition of $\BB_n$ have obvious
graphical interpretations. Hence, we have a surjective morphism
$\Psi_n:\BB_n(R)\rightarrow G\BB_n(R)$. 
It is important to note that $G\BB_n$ is, in contrast to, say, the
Temperley-Lieb algebra, not defined by giving a linear basis.
It is, rather, an algebra defined by generators and relations where not all
relations are stated explicitly.
The existence of $\Psi_n$ tells us that $2^n(2n-1)!!$ is an upper bound
for the dimension of $G\BB_n$ as well. Furthermore,
versions of propositions \ref{wortsatz} and \ref{yswort}
hold as well for this algebra.

The classical limit of a tangle algebra is defined by forgetting over
and under crossings. In our situation this should only be applied to the
crossings $X^{(G)}_i$. Then, one has $X_i^{(G)}=X_i^{(G)-1}$ and 
we demand that we have ${Y^{(G)}}^2=1$ in the limit as well.
Thus $\Psi_n(Y'_i)=\Psi_n(Y_i)$ in the limit. This shows that in the limit 
$Y^{(G)}$  behaves natural with respect to crossings and may 
therefore be represented by a dot on the arc.
Relation (\ref{lem2c}) together with $Y_i^{(G)}={Y_i^{(G)}}'=Y_i^{(G)-1}$ 
shows that in the classical limit one has 
$\Psi_n(e_iY_i)=\Psi_n(e_iY_{i+1})$.

The classical limit may be obtained by specializing the parameters
of the algebra. It is given by
\begin{eqnarray}
\BB^c_n&:=&\BB_n(R_0)\otimes _{R_0} R_c\\
R_c&:=& R_0/(\lambda-1,q-1,q_1)
\end{eqnarray}

It is obvious that $\Psi_n(\BB^c_n)$  is an algebra of 
dotted Brauer graphs. Each arc may have none or one dot on it.
Upon multiplication the number of dots is reduced modulo 2 and
a dotted cycle  is eliminated at the expense of a factor $A$. 
At the moment, however, we don't know if one obtains the full
$2^n(2n-1)!!$ dimensional dotted Brauer algebra since it 
may be that $\BB_n$ is too small.

\unitlength1mm

\section{Conditional Expectation and trace on $\BB_n$}
\label{trsec}

The graphical Interpretations suggests that a Markov trace should
exist on $\BB_n$. It will be defined as iteration of the conditional
expectation which, graphically speaking, closes the last strand.

We will need the following assumption:

\begin{hypo}\label{wenzlhypo}
The inclusion $i:\BB_n\rightarrow\BB_{n+2},
a\mapsto x^{-1}ae_{n+1}$ is injective.
 \end{hypo}

\begin{lemma}
This hypothesis is valid for $G\BB_n(R)$, that is
the morphism  $i^{(G)}:\BB_n^{(G)}\rightarrow\BB_{n+2}^{(G)},
a\mapsto x^{-1}ae^{(G)}_{n+1}$ is injective. \label{gwenzlincl}
\end{lemma}
\begin{bew}
Assume that  $a$ lies in the kernel of $i^{(G)}$.
Now, we deform the $n$-th strand of $a$ above and below of $a$
in the way indicated in figure \ref{inklu}.
Thus we have an isotopy to a graph that looks locally like $ae_{n+1}$.
So  $ae_{n+1}=0$ implies   $a=0$.
\end{bew}

   \unitlength1mm
 \begin{figure}[ht]
\begin{picture}(120,30)

\linethickness{0.1mm}

\put(10,10){\line(1,0){10}}
\put(10,10){\line(0,1){10}}
\put(20,10){\line(0,1){10}}
\put(10,20){\line(1,0){10}}
\put(12,12){\mbox{$a$}}

\put(11,20){\line(0,1){8}}
\put(13,20){\line(0,1){8}}
\put(15,20){\line(0,1){8}}
\put(17,20){\line(0,1){8}}
\put(19,20){\line(0,1){8}}

\put(11,2){\line(0,1){8}}
\put(13,2){\line(0,1){8}}
\put(15,2){\line(0,1){8}}
\put(17,2){\line(0,1){8}}
\put(19,2){\line(0,1){8}}

\put(30,12){\mbox{$=$}}

\put(40,10){\line(1,0){10}}
\put(40,10){\line(0,1){10}}
\put(50,10){\line(0,1){10}}
\put(40,20){\line(1,0){10}}
\put(42,12){\mbox{$a$}}

\put(41,20){\line(0,1){8}}
\put(43,20){\line(0,1){8}}
\put(45,20){\line(0,1){8}}
\put(47,20){\line(0,1){8}}
\put(49,20){\line(0,1){4}}

\put(41,2){\line(0,1){8}}
\put(43,2){\line(0,1){8}}
\put(45,2){\line(0,1){8}}
\put(47,2){\line(0,1){8}}
\put(49,6){\line(0,1){4}}

\put(52,6){\oval(6,6)[b]}
\put(52,24){\oval(6,6)[t]}

\put(55,10){\line(0,-1){4}}
\put(55,20){\line(0,1){4}}

\put(58,10){\oval(6,6)[t]}
\put(58,20){\oval(6,6)[b]}

\put(61,10){\line(0,-1){10}}
\put(61,20){\line(0,1){10}}

\end{picture}
\caption{\label{inklu}  }
\end{figure}

Consider 
$w=w_1\gamma w_2\in\BB_{n+1}$ with $w_i\in\BB_n,\gamma\in\Gamma_{n+1}$. 
Then we have
$e_{n+1}we_{n+1}=w_1e_{n+1}\gamma e_{n+1}w_2=sw_1w_2e_{n+1}$,
with a factor $s$ which assumes the values $s=x,1,\lambda^{-1},A$ if
$\gamma=1,e_n,X_n,Y_{n+1}$.
Thanks to hypothesis \ref{wenzlhypo} we can give the following definition
of the conditional expectation.
\begin{de}
$\epsilon_n:\BB_{n+1}\rightarrow\BB_n$ \deu{sei definiert durch}
\eng{is defined by}
$e_{n+1}ae_{n+1}=:x\epsilon_n(a)e_{n+1}$.
\end{de}
Obviously,
$\epsilon_n(w_1aw_2)=w_1\epsilon_n(a)w_2$\deu{, falls}\eng{ if} $w_i\in\BB_n$.
Furthermore, it follows from (\ref{lem1l}) 
that $e_{n+1}=e_{n+1}e_ne_{n+1}=x\epsilon_n(e_n)e_{n+1}$ thus
$\epsilon_n(e_n)=x^{-1}$. 
\deu{Analog findet man aus}\eng{Similarly one derives from} (\ref{lem1j})
\eng{the relation}
$e_{n+1}=\lambda^\pm e_{n+1}X^\pm_ne_{n+1}=
\lambda^\pm x\epsilon_n(X^\pm_n)e_{n+1}$ \deu{also}\eng{thus}
$\epsilon_n(X_n^\pm)=x^{-1}\lambda^\mp$
\deu{und aus}\eng{and from} (\ref{lem2ds}) \deu{folgt}\eng{it follows that}
$e_{n+1}=A^{-1}e_{n+1}Y_{n+1}e_{n+1}=
A^{-1}x\epsilon_n(Y_{n+1})e_{n+1}$\deu{, also}\eng{thus}
$\epsilon_n(Y_{n+1})=Ax^{-1}$. 

\eng{The itarated application of the conditional expectation yields
a map to the ground ring that will turn out to be a trace.}

\begin{de}
$\tr(a):=\tr(\epsilon_{n-1}(a)), tr(1):=1$
\end{de}
\begin{lemma}$\tr(e_n)=\epsilon_n(e_n)=x^{-1},\quad
              \tr(X_n^\pm)=\epsilon_n(X_n^\pm)=x^{-1}\lambda^\mp,\quad
              \tr(Y_{n+1})=\epsilon_n(Y_{n+1})=Ax^{-1}$
\end{lemma}  
\begin{lemma} \deu{F\"ur }\eng{$\forall$}$w_1,w_2\in\BB_n,\gamma\in\Gamma_{n+1}$ 
\deu{gilt}\eng{we have}
$\tr(w_1\gamma w_2)=\tr(\gamma)\tr(w_1w_2)$ \eng{and}
$\epsilon_n(w_1\gamma w_2)=\tr(\gamma)w_1w_2$.
\end{lemma}
\begin{bew} \deu{Die erste Aussage folgt aus der zweiten.}
\eng{The first statement is a consequence of the second which is established
in the following calculation.}
$x\epsilon_n(w_1\gamma w_2)e_{n+1}=e_{n+1}w_1\gamma w_2e_{n+1}=
w_1e_{n+1}\gamma e_{n+1}w_2=w_1x\epsilon_n(\gamma)e_{n+1}w_2=
w_1w_2x\epsilon_n(\gamma)e_{n+1}$. \deu{Daraus folgt die zweite Behauptung.	}
\end{bew}

\begin{lemma} \deu{F\"ur alle}\eng{For all} $a\in\BB_n$ \deu{gilt:}
\eng{the following equations hold.}
\begin{eqnarray}
\epsilon_n(X_n^{-1}aY'_{n+1})&=&
\epsilon_n(X_n^{-1}Y'_{n+1})a=x^{-1}\lambda^{-1}Y'_na\\
\epsilon_n(X_nY'_{n+1})&=&
x^{-1}\lambda^{-1}Y'_n+\delta Ax^{-1}-\delta\lambda x^{-1} {Y'_n}^{-1}
\end{eqnarray}
\end{lemma}
\begin{bew}
\begin{eqnarray*}
\epsilon_n(X_n^{-1}aY'_{n+1})&=&\epsilon_n(X_n^{-1}Y'_{n+1}a)\\
&=&\epsilon_n(X_n^{-1}Y'_{n+1})a=\epsilon_n(Y'_{n}X_n)a=
Y'_n\epsilon(X_n)a
=x^{-1}\lambda^{-1}Y'_na\\
\epsilon_n(X_nY'_{n+1})&=&\epsilon_n(X_n^2Y'_nX_n)=
\epsilon_n(Y'_nX_n)+\delta\epsilon_n(X_nY'_nX_n)-
\delta\lambda\epsilon_n(e_nY'_nX_n)\\
&=&Y'_n\epsilon_n(X_n)+\delta\epsilon_n(Y'_{n+1})-
\delta\lambda\epsilon_n(e_n{Y'_n}^{-1})\\
&=&x^{-1}\lambda^{-1}Y'_n+\delta Ax^{-1}-
\delta\lambda x^{-1} {Y'_n}^{-1}
\end{eqnarray*}
\end{bew}

\begin{lemma}
$\forall a\in\BB_n\quad
\epsilon_n(X_n^{-1}aX_n)=\epsilon_n(X_naX_n^{-1})=\epsilon_n(e_nae_n)=\epsilon_{n-1}(a)$
\end{lemma}
\begin{bew}
By linearity and proposition \ref{wortsatz} it is enough to show: 
\begin{eqnarray}
e_{n+1}(X_n^{-1}\gamma X_n)e_{n+1}&=&e_{n+1}(X_n\gamma X_n^{-1})e_{n+1}
\nonumber\\
&=&e_{n+1}(e_n\gamma e_n)e_{n+1}=x\tr(\gamma)e_{n+1}\nonumber
\end{eqnarray}
This is obviously true for $\gamma=1$.  For $\gamma=e_{n-1}$
one obtains
\begin{eqnarray*}
\lefteqn{e_{n+1}(X_n^{-1}e_{n-1} X_n)e_{n+1}=e_{n+1}(X_ne_{n-1} X_n^{-1})e_{n+1}=
e_{n+1}(e_ne_{n-1} e_n)e_{n+1}=xx^{-1}e_{n+1}}\\
&\Leftrightarrow&
e_{n+1}(X_{n-1}e_{n} X_{n-1}^{-1})e_{n+1}=e_{n+1}(X_{n-1}^{-1}e_{n} X_{n-1})e_{n+1}=
e_{n+1}e_ne_{n+1}=e_{n+1}\end{eqnarray*}
\deu{ und das ist unter Beachtung von (\ref{lem1i})
offensichtlich richtig.}\eng{This is true by (\ref{lem1i}).}

\deu{Im Fall}\eng{If} $\gamma=Y_n$ \deu{entsteht}\eng{one has}
\begin{eqnarray*}
\lefteqn{e_{n+1}(X_n^{-1}Y_n X_n)e_{n+1}=e_{n+1}(X_nY_n X_n^{-1})e_{n+1}=
e_{n+1}(e_nY_n e_n)e_{n+1}=x\tr(Y_n)e_{n+1}}\\
&\Leftrightarrow
e_{n+1}(X_n^{-1}Y_n X_n)e_{n+1}=e_{n+1}Y_{n+1}e_{n+1}=
e_{n+1}(e_nY_n e_n)e_{n+1}=Ae_{n+1}
\end{eqnarray*}

That this is true may be seen by transforming the first expression
\begin{eqnarray*}
\lefteqn{e_{n+1}X_n^{-1}Y_nX_ne_{n+1}=e_{n+1}e_nX_{n+1}Y_nX_ne_{n+1}=
e_{n+1}e_nY_nX_{n+1}X_ne_{n+1}=}\\
&=&
e_{n+1}e_nY_ne_nX_{n+1}X_n=Ae_{n+1}e_nX_{n+1}X_n=Ae_{n+1}
\end{eqnarray*}

The last case is $\gamma=X_{n-1}$. 
\begin{eqnarray*}
\lefteqn{e_{n+1}(X_n^{-1}X_{n-1} X_n)e_{n+1}=
e_{n+1}(X_nX_{n-1} X_n^{-1})e_{n+1}=}\\
&&=e_{n+1}(e_nX_{n-1} e_n)e_{n+1}=x\tr(X_{n-1})e_{n+1}\\
&\Leftrightarrow&
e_{n+1}(X_{n-1}X_{n} X_{n-1}^{-1})e_{n+1}=
e_{n+1}(X_{n-1}^{-1}X_{n} X_{n-1})e_{n+1}=\\&&\quad=
e_{n+1}(\lambda^{-1}e_n)e_{n+1}=\lambda^{-1}e_{n+1}\\
&\Leftrightarrow&
X_{n-1}e_{n+1}X_ne_{n+1}X_{n-1}^{-1}=X_{n-1}^{-1}e_{n+1}X_ne_{n+1}X_{n-1}=
\lambda^{-1}e_{n+1}=\lambda^{-1}e_{n+1}\\
&\Leftrightarrow&
X_{n-1}\lambda^{-1}e_{n+1} X_{n-1}^{-1}=
X_{n-1}^{-1}\lambda^{-1}e_{n+1}X_{n-1}=
\lambda^{-1}e_{n+1}
\end{eqnarray*}
\end{bew}

Now we show that $\tr$ is really a trace, i.e. $\tr(ab)=\tr(ba)$.

\begin{lemma} Assume $I_{n+1}$  to be semisimple and $\tr$ 
\deu{auf}\eng{to be a trace on} $\BB_n$
\deu{ eine Spur ist, dann ist}\eng{Then} $\tr$ 
\deu{auch auf}\eng{is a trace on} $\BB_{n+1}$\deu{ eine Spur}. \label{spurlem}
\end{lemma}
\begin{bew}
\deu{Es reicht zu zeigen, da\3}\eng{It suffices to show that}
$\tr(uv)=\tr(vu)\forall u,v\in\BB_{n+1}$. 
\deu{Falls einer der Faktoren, etwa $u$, schon
in $\BB_n$ liegt, ergibt sich das aus folgender Rechnung:}
\eng{If one of the factors, $u$ say, is actually in $\BB_n$
this follows from a simple calculation:}
$\tr(uv)=\tr(\epsilon_n(uv))=\tr(u\epsilon_n(v))=
\tr(\epsilon_n(v)u)=\tr(\epsilon_n(vu))=\tr(vu)$.

\deu{Nach Satz}\eng{Using proposition} \ref{wortsatz} 
\deu{und seinen Korollaren lassen sich
allgemeine Elemente $u,v\in\BB_{n+1}$ schreiben als}
\eng{one can write arbitrary elements $u,v\in\BB_{n+1}$ in the form}
\begin{eqnarray}
u&=&u_1+u_2Y'_{n+1}+u_3e_nu_4+u_5X_nu_6\\
v&=&v_1+v_2Y'_{n+1}+v_3e_nv_4+v_5X_n^{-1}v_6
\end{eqnarray}
\deu{Wegen der Linearit\"at von $\tr$ reicht es also aus, alle
Kombinationen der Summanden zu untersuchen. Von den 16 Kombinationen
wurden oben schon die abgehandelt, in denen einer der Faktoren in $\BB_n$
liegt. Damit bleiben neun F\"alle. Wir untersuchen zuerst die diagonalen
Kombinationen. Dabei schreiben wir stets $a$ (resp. $b$) f\"ur einen
Summanden in $u$ (resp. $v$) und bennen die $u_i,v_i$ nach Belieben um.}
\eng{Since $\tr$ is linear it suffices to proof the proposition
for all combinations. We have already dealt with the cases $u\in\BB_n$
or $v\in\BB_n$ so only nine cases remain. We investigate symmetric 
combinations first and write $a$ (resp. $b$) for one of the 
summands of $u$ (resp. $v$) and rename the $u_i,v_i$ in a handy way.}

\deu{Der erste Fall sei}\eng{First case:} 
$a=a_1e_na_2,b=b_1e_nb_2,a_i,b_i\in\BB_n$.
\begin{eqnarray*}
\tr(ab)&=&\tr(\epsilon_n(a_1e_na_2b_1e_nb_2))=\tr(a_1\epsilon_n(e_na_2b_1e_n)b_2)\\
&=&\tr(a_1\epsilon_{n-1}(a_2b_1)b_2)=\tr(b_2a_1\epsilon_{n-1}(a_2b_1))\\
&=&\tr(\epsilon_{n-1}(b_2a_1)\epsilon_{n-1}(a_2b_1))=
  \tr(\epsilon_{n-1}(a_2b_1)\epsilon_{n-1}(b_2a_1))\\
&=&\tr(a_2b_1\epsilon_{n-1}(b_2a_1))=\tr(b_1\epsilon_{n-1}(b_2a_1)a_2)\\
&=&\tr(b_1\epsilon_{n}(e_nb_2a_1e_n)a_2)=
\tr(\epsilon_n(b_1e_nb_2a_1e_na_2))=\tr(ba)
\end{eqnarray*}

\deu{Sei nun}\eng{Second case:}  $a=a_1X_{n}a_2,b=b_1X_{n}^{-1}b_2$\deu{ sein.}
\begin{eqnarray*}
\tr(ab)&=&\tr(a_1X_na_2b_1X_n^{-1}b_2)=
             \tr(a_1\epsilon_n(X_2a_2b_1X_n^{-1})b_2)\\
             &=&\tr(a_1\epsilon_{n-1}(a_2b_1)b_2)=\tr(a_1a_2b_1b_2)\\
             &=&\tr(b_1b_2a_1a_2)=\tr(ba)
\end{eqnarray*}
\deu{Der letzte Diagonale Fall ist}\eng{Third case:} 
$a=a_1Y'_{n+1},b=b_1Y'_{n+1}$.
\begin{eqnarray*}
\tr(ab)&=&\tr(a_1Y'_{n+1}b_1Y'_{n+1})=
\tr(a_1\epsilon_n(Y'_{n+1}b_1Y'_{n+1}))\\
&&\tr(a_1\epsilon_n(Y'_{n+1}Y'_{n+1}b_1))=
\tr(a_1\epsilon_n({Y'_{n+1}}^2)b_1)\\
&&\tr(b_1a_1\epsilon_n({Y'_{n+1}}^2))=
\tr(b_1\epsilon_n({Y'_{n+1}}^2)a_1)\\
&&\tr(a_1b_1\epsilon_n({Y'_{n+1}}^2))=\tr(ba)
\end{eqnarray*}
Here we used the fact that
$\epsilon_n({Y'_{n+1}}^2)$ commutes with $a_1$ since
for all $c\in\BB_n$ one has
\begin{eqnarray*}
c\epsilon_n({Y'_{n+1}}^2)e_{n+1}&=&cx^{-1}e_{n+1}{Y'_{n+1}}^2e_{n+1}
=x^{-1}e_{n+1}{Y'_{n+1}}^2e_{n+1}c\\
&=&\epsilon_n({Y'_{n+1}}^2)e_{n+1}c
\end{eqnarray*}

Fourth case:
$a=a_1Y'_{n+1},b=a_3X_n^{-1}a_4$
\begin{eqnarray*}
\tr(ab)&=&\tr(a_1\epsilon_n(Y'_{n+1}a_3X_n^{-1})a_4)=
          \tr(a_1a_3\epsilon_n(Y'_{n+1}X_n^{-1})a_4)\\
&=&x^{-1}\lambda^{-1}\tr(a_1a_3Y'_na_4)=x^{-1}\lambda^{-1}\tr(a_3Y'_na_4a_1)\\
&=&\tr(a_3\epsilon_n(X_n^{-1}Y'_{n+1})a_4a_1)=
\tr(\epsilon_n(a_3X_n^{-1}a_4a_1Y'_{n+1}))=\tr(ba)
\end{eqnarray*}

\deu{N\"achster Fall}\eng{Sixth case}: $a=a_1X_na_2,b=a_3Y'_{n+1}$.
\begin{eqnarray*}
\tr(ab)&=&\tr(a_1\epsilon_n(X_na_2a_3Y'_{n+1}))=
 \tr(a_1\epsilon_n(X_nY'_{n+1})a_2a_3)\\
&=&x^{-1}\lambda^{-1}\tr(a_1Y'_{n}a_2a_3)+\delta A x^{-1}\tr(a_1a_2a_3)-
\delta\lambda x^{-1}\tr(a_1{Y'_n}^{-1}a_2a_3)\\
&=&x^{-1}\lambda^{-1}\tr(a_3a_1Y'_{n}a_2)+\delta A x^{-1}\tr(a_3a_1a_2)-
\delta\lambda x^{-1}\tr(a_3a_1{Y'_n}^{-1}a_2)\\
&=&\tr(a_3a_1\epsilon_n(Y'_{n+1}X_n)a_2)=\tr(ba)
\end{eqnarray*}

\deu{Sei nun}\eng{Seventh case:} $a=a_1e_na_2,b=a_3Y'_{n+1}$.
\begin{eqnarray*}
\tr(ab)&=&\tr(a_1\epsilon_n(e_na_2a_3Y'_{n+1}))=
  \tr(a_1\epsilon_n(e_nY'_{n+1})a_2a_3)\\
&=&\lambda\tr(a_1\epsilon_n(e_n{Y'_n}^{-1})a_2a_3)=
\lambda\tr(a_1\epsilon_n(e_n){Y'_n}^{-1}a_2a_3)\\
&=&\lambda x^{-1}\tr(a_1{Y'_n}^{-1}a_2a_3)=
   \lambda x^{-1}\tr(a_2a_3a_1{Y'_n}^{-1})\\
&=&\lambda\tr(a_2a_3a_1\epsilon_n(e_n){Y'_n}^{-1})=
\tr(a_2a_3a_1\epsilon_n(Y'_{n+1}e_n))\\
&=&\tr(\epsilon_n(a_3Y'_{n+1}a_1e_n)a_2)=\tr(ba)
\end{eqnarray*}
\deu{Damit ist auch der Fall}\eng{The case} $b=a_1e_na_2,a=a_3Y'_{n+1}$
\deu{abgehandelt}\eng{is similar}.
\deu{Es bleiben jetzt nur noch die nicht diagonale F\"alle, in 
denen ein Faktor $e_n$ enth\"alt. Diese treten aber defacto
gar nicht auf.}
\eng{The only remaining cases are nonsymmetric with one occurrence of $e_n$.}
\deu{Denn wegen der vorausgesetzten 
Halbeinfachheit des Ideals}
\eng{Since we assume} $I_{n+1}$ 
\deu{gibt es eine zentrale Idempotente}\eng{to be semisimple there is an idempotent}
$z\in\BB_{n+1}$ \deu{mit}\eng{such that}
$z\BB_{n+1}\iso I_{n+1}$. \deu{Sei nun}\eng{Now assume that $a$ contains
$e_n$, hence} $a\in I_{n+1}$\deu{, also}\eng{ i.e.} $a=az$. 
\deu{Dann ist}\eng{Then we have}
$ab=azb=a(zb)$\deu{. Also kann ObdA auch}\eng{, which shows that we might 
as well assume} $b\in I_{n+1}$\eng{.} 
\deu{angenommen werden (ansonsten ist die Spur $0$).}
\deu{Aber}\eng{But} $a,b\in I_{n+1}$ 
\deu{impliziert, da\3 es sich um Linearkominationen der Form}
\eng{implies that $a,b$ are linear combinations of the form}
$a=\sum_ia_ie_na'_i,b=\sum_ib_ie_nb'_i$ \deu{mit}\eng{with}
$a_i,a'_i,b_i,b'_i\in\BB_n$\deu{ handelt}. \deu{Terme
dieser Art wurden aber schon behandelt.}\eng{Thus we are back in a 
case that was already treated.}
\end{bew}

\eng{\section{The structure theorem}}  \label{mainsec}

We only need a few definitions on Young diagrams before we can state
the structure theorem for $\BB_n$.

A Young diagram $\lambda$ of size $n$ 
is a partition of the natural number $n$.
$\lambda=(\lambda_1,\ldots,\lambda_k),\sum_i\lambda_i=n,
\lambda_i\geq\lambda_{i+1}$. 
In the following we use ordered pairs of Young diagrams
(cf. \cite{ariki}). The size of a pair of Young diagrams is the sum of 
sizes of its components.
Let $\widehat{\Gamma}_n$ be the set of all pairs of Young diagrams of sizes 
$n,n-2,\ldots$.

\begin{satz}
The following statements hold for the algebra
$\BB_n(K_0)$ over the quotient field $K_0$.
\begin{enumerate}
\item $\BB_n$ is isomorphic to $G\BB_n$ and it is semisimple.
The simple components are indexed by
$\widehat{\Gamma}_n$. \label{decomp}
\begin{equation}\BB_n=\bigoplus_{(\mu,\lambda)\in\widehat{\Gamma}_n}
 \BB_{n,(\mu,\lambda)}\end{equation}
\item The Bratteli rule for restrictions of modules: A simple
$\BB_{n,(\nu,\rho)}$ module 
$V_{(\nu,\rho)},(\nu,\rho)\in\widehat{\Gamma}_n$ 
 decomposes into $\BB_{n-1}$ 
modules such that the $\BB_{n-1}$ module  
$(\mu,\lambda)\in\widehat{\Gamma}_{n-1}$ occurs iff
$(\mu,\lambda)$ may be obtained from $(\nu,\rho)$ 
by adding or removing a box.
\item $\tr$ is a faithful trace. 
To every pair of Young diagrams 
$(\mu,\lambda)\in\widehat{\Gamma}_n$ 
there is a minimal idempotent
$p_{(\mu,\lambda)}$ and a non vanishing, rational function	 
$Q_{(\mu,\lambda)}$ which does not depend on $n$ and satisfies
$\tr(p_{(\mu,\lambda)})=Q_{(\mu,\lambda)}/x^n$.
\end{enumerate}
\end{satz}

  \unitlength1mm
 \begin{figure}[ht]
\begin{picture}(130,30)

\put(10,20){\mbox{$\BB_0$}}
\put(10,10){\mbox{$\BB_1$}}
\put(10,0){\mbox{$\BB_2$}}

\put(70,20){\mbox{$(\cdot,\cdot)$}}

\put(60,10){\mbox{$(\Box,\cdot)$}}
\put(80,10){\mbox{{\small $(\cdot,\Box)$}}}
\put(73,18){\line(1,-1){5}}
\put(73,18){\line(-1,-1){5}}

\put(65,0){\mbox{{\small $(\cdot,\cdot)$}}}
\put(75,0){\mbox{{\small $(\Box,\Box)$}}}
\put(35,0){\mbox{{\small $(\Box\Box,\cdot)$}}}
\put(50,0){\mbox{{\small $({\Box\atop\Box},\cdot)$}}}
\put(90,0){\mbox{{\small $(\cdot,{\Box\atop\Box})$}}}
\put(105,0){\mbox{{\small $(\cdot,\Box\Box)$}}}

\put(63,8){\line(-4,-1){20}}
\put(63,8){\line(1,-1){5}}
\put(63,8){\line(-1,-1){5}}
\put(63,8){\line(3,-1){15}}
\put(85,8){\line(-3,-1){15}}
\put(85,8){\line(1,-1){5}}
\put(85,8){\line(-1,-1){5}}
\put(85,8){\line(4,-1){20}}

\linethickness{0.2mm}

\end{picture}
\caption{\label{brat} 
\deu{Das Bratteli-Diagram von $\BB_n$}
\eng{The Bratteli digram of $\BB_n$}}

\end{figure}
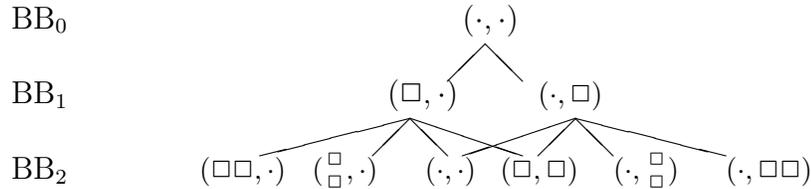

For the proof of the structure theorem we need  some  facts from
Jones-Wenzl theory of
inclusions of finite dimensional semisimple algebras.

Let  $A\subset B\subset C$ 
be a unital imbedding of finite dimensional semisimple algebras and let
$\tr$ be a trace on  $A, B$ 
that is compatible with the inclusion. 
\eng{The associated conditional expectation is denoted by} 
$\epsilon_A:B\rightarrow A,\tr(ab)=\tr(a\epsilon_A(b))$.
\eng{It is assumed that there is an idempotent $e\in C$ such that }
$e^2=e,ebe=e\epsilon_A(b)\forall b\in B$\deu{, }\eng{ and } 
$\varphi:A\rightarrow C,a\mapsto ae \mbox{ \deu{injektiv}\eng{is injective}}$.

\eng{Such a situation can be realized starting from an inclusion pair }
$A\subset B$ with a common faithful trace  $\tr$ 
and conditional expectation $\epsilon_A$. We set 
$\widehat{C}:=\{\alpha:B\rightarrow B\mid \alpha\mbox{ linear},
\alpha(ba)=\alpha(b)a\forall a\in A,b\in B\}$.
The inclusion  $B\subset\widehat{C}$ is given by
$b\mapsto\alpha_b,\alpha_b(b_1):=bb_1$. Here  $e$ is given by  
$e_A=\epsilon_A:B\rightarrow B$.
The subalgebra of $\widehat{C}$ generated by  $B$ and  $e_A$ 
is denoted by $<B,e_A>$.
\eng{For this setup Wenzl has obtained the following results} 
\cite[Theorem 1.1]{we2}
\begin{enumerate}
\item $<B,e_A>\iso{\rm End}_A(B)$\label{towera}
\item \deu{Die einfachen Komponenten von}\eng{The simple components of} 
$A$ \deu{und}\eng{and} $<B,e_A>$ \deu{sind bijektiv
zuordenbar}\eng{are in 1-1 correspondence}. 
The inclusion matrices of
$A\subset B\subset<B,e_A>$ are relatively transposed. 
If $p$ is a minmal idempotent in $A$ then $pe_A$ is minimal idempotent
in $<B,e_A>$. 
\label{towerb}
\item $<B,e_A>\iso Be_AB$ \label{towerc}
\item $<B,e>\iso<B,e_A>\oplus\widetilde{B}$ where $\widetilde{B}$
\eng{is a subalgebra of} $B$\deu{ ist}. \label{towerd}
\item \eng{\ref{towerd} implies that the ideal genarated by}  $e$ in $C$ 
is isomorphic to $<B,e_A>$.
\label{towere}
\end{enumerate}

We now give the proof of the main theorem.

\begin{bew}
$\BB_0$ is simply the ground ring. 
\eng{Thus the proposition is true with}
$\tr(p_{(\cdot,\cdot)})=\tr(1)=Q_{(\cdot,\cdot)}/x^0,
Q_{(\cdot,\cdot)}=1$. 
The algebra $\BB_1$ is twodimensional and has a basis $\{1,Y\}$. 

Assume the proposition is shown by induction for $\BB_n$.

By the induction assumptions we have $\BB_n=G\BB_n$. Using this we show that
the inclusion $i:\BB_n\rightarrow\BB_{n+2}$ of section \ref{trsec}
is injective. Assume we have $i(a)=0$, then 
$0=\Psi_{n+2}(i(a))=i^{(G)}(a)$ and the claim follows from
injectivity of $i^{(G)}$. 

We apply the Jones-Wenzl theory to the following situation:
$A=\BB_{n-1},B=\BB_n,C=\BB_{n+1},e=x^{-1}e_n,\epsilon_A=\epsilon_{n-1}$.
This is possible because $A,B$ are semisimple algebras with a
faithful trace by induction assumption.
All properties needed for $e$ have already been established.
Statement \ref{towera} of Jones-Wenzl theory asserts the semisimplicity of
${\rm End}_A(B)\iso<B,e_A>$  
which is by \ref{towere} the ideal generated by $e$.
Thus $I_{n+1}$ is semisimple.
The quotient algebra by $\BB_{n+1}/I_{n+1}$ is the Hecke algebra $\HB_{n+1}$
and is semisimple according to \cite{ariki}.
Now, in general if $A$ is a finite dimensional algebra over some field
with a semsisimple ideal $I$ such
that $A/I$ is semisimple as well then $A$ is semisimple itself:
The map $A\rightarrow A/I$ maps the radical ${\rm Rad}(A)$ into the
radical of $A/I$ which is tivial, hence ${\rm Rad}(A)\subset I$ and 
thus ${\rm Rad}(A)=I\cap{\rm Rad}(A)\subset{\rm Rad}(I)=\{0\}$. 
For finite dimensional algebras over a field vanishing of the radical
is equivalent to semisimplicity.

Thus $\BB_{n+1}$ is semisimple and is a direct sum 
$\BB_{n+1}=I_{n+1}\oplus\BB_{n+1}/I_{n+1}$.
Now, the same reasoning can be applied to the the algebra $G\BB_n$.
In this case the quotient  $G\BB_{n+1}/I^{(G)}_{n+1}$ arises. 
Imposing the relation $e^{(G)}_i=0$ obviously annihilates all
tangles that are not ribbon braids of B-type. But then standard knowledge
about the graphical interpretation of Hecke algebras 
shows that $\HB_{n+1}=G\BB_{n+1}/I^{(G)}_{n+1}$ as well.
Jones-Wenzl theory then implies $G\BB_{n+1}=\BB_{n+1}$.

Statement \ref{towerb} asserts that the simple components of 
$I_{n+1}$ are indexed by $\widehat{\Gamma}_{n-1}$.
The simple components  of $\HB_{n+1}$ are indexed by pairs of
Young diagrams of size $n+1$ (see \cite{ariki}). 
This completes the proof of point \ref{decomp} of the theorem.

The inclusion matrix for the part $I_{n+1}$ is the transpose
of the inclusion matrix of $\BB_{n-1}\subset\BB_n$. For the part $\HB_{n+1}$
the Bratteli rule follow from \cite{ariki}.

The results proven sofar and lemma \ref{spurlem} imply that $\tr$
is a trace. To show its  faithfulness one has to show that the $Q$
functions don't vanish.
If $p_{(\mu,\lambda)}\in\BB_{n-1}$ is a minimal idempotent in
$\BB_{n-1,(\mu,\lambda)}$ then $x^{-1}p_{(\mu,\lambda)}e_n$ 
is a minimal idempotent in $\BB_{n+1}$.
The trace of this idempotent is
$\tr(x^{-1}p_{(\mu,\lambda)}e_n)=x^{-2}tr(p_{(\mu,\lambda)})=
Q_{(\mu,\lambda)}/x^{n-1+2}$.
Obviously, this is nonvanishing (using the induction assumption). 
The idempotents of this kind are those of $I_{n+1}$.
For the other idempotents (which are those of $\BB_{n+1}/I_{n+1}$)
the function $Q$ is defined by
$\tr(p_{(\mu,\lambda)})=Q_{(\mu,\lambda)}/x^n$.

Now, we have two possibilities to establish faithfulness of the trace.
One way is to note that  $\tr$ restricted to  $\HB_{n+1}=\BB_{n+1}/I_{n+1}$ 
is the Markov trace of the Hecke algebra which is known to be nondegenerate.
The second possibility is to use the classical limit.
A minimal idempotent $p_{(\lambda,\mu)}$ of $\BB_n$ yields
an idempotent in the classical limit. 
On this algebra the trace in known to be nondegenerate \cite{reich}. 
Thus the function $Q_{(\lambda,\mu)}$ has a non vanishing limit.  
\end{bew}

In the rest of this section we sketch a second proof
of the semisimplicity of $\BB_n(K_0)$.
It is based on a different approach to the Markov trace
which is based on a different realization of process
that may graphically interpreted as closing tangles.

We start with some definitions.
\begin{eqnarray}
X(i,j)&:=&X_iX_{i+1}\cdots X_j\\
X^{-1}(i,j)&:=&X_i^{-1}X_{i+1}^{-1}\cdots X_j^{-1}\\
E(i,j)&:=&e_ie_{i+2}\cdots e_j\\
H_1&:=&e_1\\
H_{n+1}&:=&e_{n+1}X(n+2,2n+1)X(n+1,2n)H_n
\end{eqnarray}

\unitlength1mm
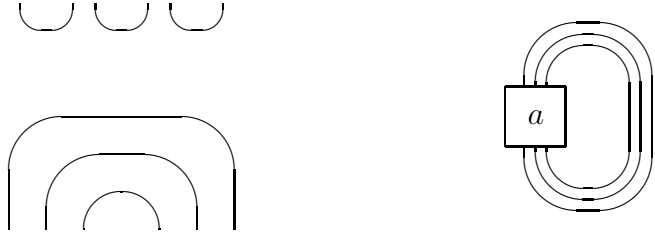
\begin{figure}[ht]
\begin{picture}(120,30)

\put(20,0){\oval(30,30)[t]}
\put(20,0){\oval(20,20)[t]}
\put(20,0){\oval(10,10)[t]}

\put(20,30){\oval(7,7)[b]}
\put(10,30){\oval(7,7)[b]}
\put(30,30){\oval(7,7)[b]}

\put(71,11){\line(0,1){8}}
\put(79,11){\line(0,1){8}}
\put(71,11){\line(1,0){8}}
\put(71,19){\line(1,0){8}}
\put(74,14){\mbox{$a$}}

\put(82,19){\oval(11,11)[t]}
\put(82,19){\oval(14,14)[t]}
\put(82,19){\oval(17,17)[t]}

\put(82,11){\oval(11,11)[b]}
\put(82,11){\oval(14,14)[b]}
\put(82,11){\oval(17,17)[b]}
\put(87.5,11){\line(0,1){8}}
\put(89,11){\line(0,1){8}}
\put(90.5,11){\line(0,1){8}}

\end{picture}
\caption{\label{hbld} The graphical interpretations of $H_3$	 
(on the left) and of  $\tr(a)$ (on the right)}
\end{figure}

The following properties can be shown by
straightforward (inductive) calculations.
\begin{lemma}\begin{eqnarray}
H_n&=&E(n,n)E(n-1,n+1)\cdots E(1,2n-1)\label{hdurche}\\
H_{n+1}&=&e_{n+1}X^{-1}(n+2,2n+1)X^{-1}(n+1,2n)H_n \label{negh}\\
X_i^\pm H_n&=&X_{2n-i}^\pm H_n,\quad e_iH_n=e_{2n-i}H_n\\
e_{2n-1}&=&X(n,2n-2)^{-1}X(n+1,2n-1)^{-1}e_n\nonumber\\
  &&X(n+1,2n-1)X(n,2n-2)\\
e_n&=&X(n+1,n+k)X(n,n+k-1)e_{n+k}\nonumber\\
  &&X(n,n+k-1)^{-1}X(n+1,n+k)^{-1}\\
H_{n+1}&=&e_{n+1}X(n+2,2n)X(n+1,2n-1)X_{2n+1}^{-1}X_{2n}^{-1}H_n\\
Y^{\pm1} H_n&=&\lambda^{\pm1} {Y'_{2n}}^{\mp1}H_n\\
\overline{H_n}X_i^\pm&=&\overline{H_n}X_{2n-i}^\pm\\
\overline{H_n}e_i&=&\overline{H_n}e_{2n-i}\\
\overline{H_n}Y^{\pm1}&=&\lambda^{\pm1}
           \overline{H_n}{Y'_{2n}}^{\mp1}\\
\overline{H_n} ab H_n&=&\overline{H_n} ba H_n,\quad\forall a,b\in\BB_n\\
x^n\tr(a)E(1,2n-1)&=&\overline{H_n}aH_n\quad\forall a\in\BB_n\\
0&=&x^n(\tr(ab)-\tr(ba))E(1,2n-1)
\end{eqnarray}
\end{lemma}

Recall that $e_1$ does not vanish and has vanishing annulator ideal
in $\BB_2(R_0)$. Similarly, the same is true for $e_1^{(G)}$.
By induction using lemma \ref{gwenzlincl}
it follows that the same is true for 
 $E(1,2n-1)\in G\BB_{2n}(R_0)$. This shows that $\tr$ is 
  a trace on $G\BB_n$.

We now investigate properties of the trace in the classical
limit.
Let $a$ be a dotted Brauer graph and let $n_i(a),i=0,1$ 
be the number of cycles in its closure with $i$ dots on it.
The the trace of $a$ may easily seen to be given by
 \begin{equation}
 \tr(a)=x^{-n}x^{n_0(a)}A^{n_1(a)}
 \end{equation}

\begin{satz}
$\tr$ is nondegenerate and hence $G\BB_n(K_0)$ is semisimple.
Furthermore, $G\BB_n(K_0)=\BB_n(K_0)$.
\end{satz}
\begin{bew}
Let $S_n=\{ v_i\mid 1\leq i\leq 2^n(2n-1)!!\}$ be a set of elements 
that generate
$G\BB_n(R_0)$ and yield a basis of dotted Brauer graphs in the
classical limit. 

To prove the first statement of the proposition it is enough to show that
$0\neq{\rm det}(\tr(v_iv_j^\ast)_{i,j})\in R_0$.
We tensor this element with $R_c$ to pass
to the classical limit.
The involution $a\mapsto a^\ast$ maps graphs to their
top-down mirrored image while keeping dots.
Due to the reduction of dots modulo 2 there are no dots in 
the closure of $aa^\ast$.
Assume  $a$ has  $s$ upper (and $s$ lower) 
horizontal arcs. Then $aa^\ast$ has  $s$ cycles. 
When closing to calculate the trace another $s$ cycles
arise from the $s$ lower ans $s$ upper horizontal arcs of $a$ and $a^\ast$.
The vertical arcs of $a$ describe a permutation 
and $a^\ast$ contains the inverse permutation. 
Thus, upon closing, these vertical arcs yield another
  $n-2s$ cycles. We conclude that $\tr(aa^\ast)=1$.
Now, we specify $A=x^{-1}$ by forming a further tensor product. 
The trace will then be a Laurent polynomial in $x$.
The choice of $A$ lets dots on arcs decrease the degree of the 
trace polynomial.
Now, denote by $\beta$ an arc in  $a$ and let $b$ be another graph
which does not contain an arc that is the involutive image of $\beta$.
Investigating the cases that $\beta$ is horizontal or vertical 
one observes that the cylce in $\tr(ab)$ containing $\beta$
must contain more than  two arcs of $a$ and $b$. 
The trace of $ab$ thus is of lower degree in $x$ than the trace of
 $aa^\ast$. We conclude that
  $b=a^\ast$ is the unique graph of highest degree of $x$ in
   $\tr(ab)$. 

Using this we can establish that \[ {\rm det}(\tr(v_iv_j^\ast)_{i,j})
=x^{-nk^n(2n-1)!!}{\rm det}
((x^{n_0(v_iv_j^\ast)}x^{-n_1(v_iv_j^\ast)})_{i,j})\] 
does not vanish. The diagonal elements in this matrix are those 
of highest $x$-degree in each row.
Evaluation of the determinant thus yields only one term with highest
$x$-degree and hence the determinant cannot vanish. But then 
the original determinant of the trace on $G\BB_n(R_0)$ has to
be non zero.

The inclusion image of $S_n$ in $G\BB_n(K_0)$ generates
this algebra as a $K_0$ vector space and the determinant of the trace
is the same nonvanishing element of $R_0\subset K_0$ as before.
Existence of a nondegenerate trace on an algebra over a
a field of characteristic zero implies its semisimplicity.

A further consequence is that the dimension of $G\BB_n(K_0)$ is actually
equal to $2^n(2n-1)!!$.
The surjection $\Psi_n:G\BB_n(K_0)\rightarrow\BB_n(K_0)$ 
is thus an isomorphism.
\end{bew}

\section{Tensor representations}	 \label{repsec}

Tensor representations of $\BB_n$ were found by tom Dieck \cite{tD3}.
We review their definition and show that they can be used
to calculate the trace on $\BB_n$ as a matrix trace.
The ground field $K$ is either the function field $\CC(q)$ or $\CC$
with an element $q\in\CC$.
The construction uses the R-matrix of the quatum group
$U_q(so_N),N=2m+1,m\in\NN$. 
\eng{The $N$ dimensional defining representation operates on}
$V=\{v_i\mid i\in I\}$. 
The index set is $I=\{-N+2,-N+4,\ldots,-3,-1,0,1,3,\ldots,N-2\}$.
The permuting R-matrix is
\begin{eqnarray}
B&=&\sum_{i\neq0}(qf_{i,i}\otimes f_{i,i}+q^{-1}f_{i,-i}\otimes f_{-i,i})
+f_{0,0}\otimes f_{0,0}+\sum_{i\neq j,-j} f_{i,j}\otimes f_{j,i}+\\
&&(q-q^{-1})\left(\sum_{i<j}f_{i,i}\otimes f_{j,j}-
\sum_{j<-i}q^\frac{i+j}{2}f_{i,j}\otimes f_{-i,-j} \right)\nonumber
\end{eqnarray}
Here $f_{i,j}$ is the $N\times N$ 
matrix with a $1$ at position $(i,j)$ and $0$ elsewhere.

\deu{Aus} $E:=1-(B-B^{-1})/\delta$ is given by
\begin{equation}
E=\sum_{i,j}q^\frac{i+j}{2}f_{i,j}\otimes f_{-i,-j}
\end{equation}
This implies $E^2=xE$ with 
$x=\sum_iq^i$ and hence $\lambda=q^{1-N}$.

\eng{T. tom Dieck has found the following representaing matrix for $Y$.} 
\begin{equation}
F=-f_{0,0}+q^{-1/2}\sum_{i\neq0}f_{-i,i}+(q^{-1}-1)\sum_{i>0}f_{i,i}
\end{equation}
\eng{It satisfies} $F^2=(q^{-1}-1)F+q^{-1},
(F\otimes 1)B(F\otimes 1)B=B(F\otimes 1)B(F\otimes 1),
E=E(F\otimes1)B(F\otimes1)$.
Hence $\phi:\BB_n\rightarrow{\rm End}(V^{\otimes n}),
Y\mapsto F\otimes1\cdots\otimes 1,X_i\mapsto 1\otimes\cdots\otimes1
\otimes B\otimes1\cdots\otimes1$
defines a representation of $\BB_n$.
The parameters are $q_1=(q^{-1}-1),\lambda=q^{1-N}$.

Let $D$ be the matrix  $D_{i,i}:=q^i$ and define
$\Psi:{\rm End}(V^{\otimes n})\rightarrow K,
\Psi(a):=\Tr(a(D^{\otimes n}))/\Tr(D^{\otimes n})$.
Here ${\rm Tr}$ is the usual trace of matrices.
\begin{lemma} $\tr=\Psi\circ\phi$
\end{lemma}
\begin{bew}
Using the parameters of the tensor representation we obtain: 
\[\tr(Y)=\frac{A}{x}=\frac{q_1}{1-q_0\lambda}=
\frac{q^{-1}-1}{1-q^{-1}q^{1-N}}
=\frac{q^{-1}-1}{1-q^{-N}}\]

We now calculate $\Psi(Y)$:
\begin{eqnarray*}\Tr(D)&=&\sum_{I\ni i>0}q^i+q^0+\sum_{I\ni i<0}
=1+q^{-1}\sum_{i=1}^m(q^2)^i+q\sum_{i=1}^m(q^{-2})^i\\
&=&1+\frac{q-q^{N}}{1-q^2}+\frac{q^{-1}-q^{-N}}{1-q^{-2}}\\
\Tr(DF)&=&-1+\sum_{I\ni i>0}(q^{-1}-1)q^i=
(q^{-2}-q^{-1})\sum_{i=1}^m(q^2)^i-1\\&=&
\frac{q^{-2}-q^{-1}}{1-q^2}(q^2-q^{N+1})-1\\
\Psi(Y)&=&\Tr(DF)/\Tr(D)=
{\frac {-q^{N+1}+q^{2 N}-q^{2 N-1}+q^{N+2}}{q^{N}-q^{N+2}-q^{2 N}+q
^{2}}}
=-{\frac {q-1}{q-q^{-N+1}}}
\end{eqnarray*}
The rest of the proof coincides with the proof of \cite{we2}[Lemma 5.4]
\end{bew}

 A physical application of tensor representations  of $\BB_n$
 has been found in \cite{rhoref}.
 Two dimensional integrable systems are described by solutions of the 
spectral parameter dependent Yang-Baxter-Equation (YBE) that reads
with $R\in{\rm End}(V\otimes V)$):
\begin{equation}
R_1(t_1)R_2(t_1t_2)R_1(t_2)=R_2(t_2)R_1(t_1t_2)R_2(t_1)
\quad\forall t_1,t_2\label{ybe}
\end{equation}
If the system is restricted to a half plane an additional matrix 
$K(t)\in{\rm End}(V)$  is needed to describe reflections. Is
has to fulfill Sklyanin's reflection equation \cite{levy}:

\begin{equation}\label{re}
R(t_1/t_2)(K(t_1)\otimes1)R(t_1t_2)(K(t_2)\otimes1)=
(K(t_2)\otimes1)R(t_1t_2)(K(t_1)\otimes1)R(t_1/t_2)
\end{equation}

It is possible to obtain solutions of the YBE by Baxterization
from the A-type BMW algebra \cite{cheng}:
\begin{equation}
R_i(t)= -\delta t (t+q\lambda^{-1})+
(t-1)(t+q\lambda^{-1})X_i+\delta t(t-1)e_i
\end{equation}

Using the additional generator $Y$ of $\BB_n$ one can extend this to 
obtain solutions of the reflection equation:

\begin{satz} $K(t)=(t^2q_1(1-t^2)^{-1}+Y)f_1(t)$ 
is for arbitrary $f_1$ a solution of the reflection equation (\ref{re}).
\end{satz}

It is a remarkable fact that no similar solution 
exists for the Hecke algebra $\HB_n$.

\section{Application: Invariants of links in a solid torus.}
\label{blinks}

The Markov trace can be used to define a link invariant for 
links of B-type which are links in a solid torus. 
There is an analog of Markov's theorem for type B links found by S.
Lambrodopoulou in \cite{lamb}. It takes the same form as the usual Markov theorem,
i.e. two B-braids $\beta_1,\beta_2$ have isotopic closures
 $\hat{\beta}_1,\hat{\beta}_2$ if $\beta_1,\beta_2$
may transformed in one another by a finite sequence of moves of
the following two kinds: I Conjugation $\beta\sim\alpha\beta\alpha^{-1}$
and II $\alpha\sim\alpha\tau_n$ for $\alpha\in\ZB_n$.

This theorem implies that there exists an extension of the Kauffman polynomial 
to braids of B-type. Denote by $\pi:\ZB_n\rightarrow\BB_n$ the morphism
$\tau_i\mapsto X_i, \tau_0\mapsto Y$. Then we obtain without any further proof
an invariant of the B-type link $\hat{\beta}$
that is the closure of a B-braid $\beta_\in\ZB_n$ by the following definition:
\begin{de}
The B-type Kauffman polynomial of a B-link $\hat{\beta}$ is defined to be
\begin{equation}
L(\hat{\beta},n):=x^{n-1}\lambda^{e(\beta)}\tr(\beta)\qquad\beta\in\ZB_n
\end{equation}
$e:\ZB_n\rightarrow\ZZ$ is the exponential sum with $e(X_i)=1,e(Y)=0$.
\end{de}

\end{document}